\documentclass[draft]{svj3}                     
\smartqed  
\usepackage{graphicx}

\usepackage{amsmath}
\usepackage{amsfonts,amssymb}

\usepackage[mathcal]{eucal}
\usepackage{stmaryrd}

\usepackage[cp1251]{inputenc}   
\usepackage[english, russian]{babel}   

\newcommand{\doc}{\underline{\it Proof.}\
}
\newcommand{\bo}{\qed}
\newcommand{\rav}{\stackrel{\triangle}{=}}
\newcommand{\ph}{\varphi}
\newcommand{\la}{(}
\newcommand{\ra}{)}
\newcommand{\epsi}{\varepsilon}
\newcommand{\beg}{\varkappa}

\newcommand{\rref}[1]{$(\ref{#1})$}
\newcommand{\mm}[1]{{\mathbb{#1}}}
\newcommand{\ct}[1]{{\mathcal{#1}}}
\newcommand{\fr}[1]{{\mathfrak{#1}}}
\newcommand{\td}[1]{\widetilde{#1}}
\newcommand{\ravref}[1]{\stackrel{(\ref{#1})}{=}}
\newcommand{\leqref}[1]{\stackrel{(\ref{#1})}{\leq}}
\newcommand{\leqtref}[1]{\stackrel{\tref{#1}}{\leq}}
\newcommand{\geqtref}[1]{\stackrel{\tref{#1}}{\geq}}
\newcommand{\geqref}[1]{\stackrel{(\ref{#1})}{\geq}}

\newcommand{\tref}[1]{\textrm{$\widetilde{\textrm{\rref{#1}}}$}}
\newcommand{\rai}{\stackrel{T\to\infty}{\rightsquigarrow}}
\usepackage[active]{srcltx}
\newcommand{\ee}{
}

\journalname{Journal of Dynamical and Control Systems}
\begin{document}

\title{NECESSITY OF VANISHING SHADOW PRICE\\IN INFINITE HORIZON CONTROL PROBLEMS\thanks{}}


\titlerunning{NECESSITY OF VANISHING SHADOW PRICE}     

\author{Dmitry Khlopin 
}

\institute{D. Khlopin \at
              Institute of Mathematics and Mechanics, Ural Branch,     Russian
Academy of Sciences, 16, S.Kovalevskaja St., 620990, Yekaterinburg, Russia \\
Chair of Applied Mathematics, Institute of Mathematics and Computer Science, Ural Federal University, 4, Turgeneva St., 620083, Yekaterinburg, Russia  \\
              \email{khlopin@imm.uran.ru}           
}

\date{$\ $}

\maketitle

 \vspace{-25mm plus 0mm minus 0mm}
 This paper investigates the
necessary  optimality conditions  for uniformly overtaking
optimal
 control  on infinite horizon  in the free end case. 
 In  the  papers of    S.M.\,Aseev,  A.V.\,Kryazhimskii, V.M.\,Veliov, K.O.\,Besov
 there was suggested the boundary condition for
equations of the Pontryagin Maximum Principle. Each optimal process
 corresponds to  a  unique solution
 satisfying the boundary condition. Following A.Seierstad's idea,
 in this paper we prove a more general geometric variety of
 that boundary condition. We show that this condition is necessary
  for uniformly overtaking optimal control
   on infinite horizon in the free end case.  A number
   of assumptions under which this condition
   selects a unique Lagrange multiplier is obtained.  The results are applicable
    to general non-stationary systems and
 the optimal objective value is not necessarily
finite. Some examples are discussed.


\section*{Introduction}
\label{intro}

The Pontryagin Maximum Principle for infinite horizon problems had
already been formulated in the monograph~\cite{ppp}; the general
Maximum Principle for infinite interval was proved in \cite{Halkin},
but such Maximum Principle has no transversality condition, and in
general, selects a
too broad family of extremal trajectories.
A significant number
\cite{Halkin,kr_as,aucl,kam,mich,smirn,ssbook,sch} of such conditions
was proposed; however, as it was noted in, for example,
\cite{Halkin,mich,shell},\cite[Sect. 6]{kr_as},\cite[Example
10.2]{norv}, these conditions, as pointed out in~\cite{nnn},
 ``may frequently fail to hold (conditions
securing these properties may fail.) Even if they do hold, for
example when strong enough growth conditions hold, these condition
may fail to give any information determining the integration
constants arising when integrating the adjoint equation.''

Since the necessity of this condition does not imply its
nontriviality on solutions of the relations of the Maximum Principle,
it is reasonable to find a condition that would select a single
solution of the relations of the Maximum Principle for any optimal
control. For this purpose, \cite{norv} proposes to find $\psi^0$ such
that it is a pointwise limit of a sequence of shadow prices  equal
to  zero on certain sequence $\tau\uparrow\infty$ of times. Under
assumptions of
 \cite[Theorem 6.1]{norv}, such~$\psi^0$ is unique; in what follows,
 it will be referred to as $\tau$-vanishing shadow price.

In papers~\cite{kr_as2,kr_asD,kr_as,kr_as3}, Aseev and Kryazhimskii
proposed the analytic
  expression for the shadow prices.
  This version of the normal form of the Maximum Principle
holds with the explicitly specified shadow price. This gives a
complete set
  of necessary optimality conditions (see \cite{kr_as2,kr_asD,kr_as,kr_as3,kab,av});
moreover,  under  assumptions of
  \cite{kab,av,aucl,norv},
   the solution of this
  form of the Maximum Principle is uniquely determined by the optimal
  control.

This paper aims to merge these two approaches, to find
assumptions such that a
$\tau$-vanishing Lagrange multiplier of the Maximum Principle corresponds to
every optimal control, and to express its shadow price explicitly in the
form of an improper integral that depends only on optimal control and
trajectory.

In this paper, we consider only the problem with free right end. It
is assumed a priori that an optimal control (uniformly weakly
overtaking optimal control) exists (for discussion of existence,
refer to, for example, \cite{bald,car1,carm,car,dm}). In addition to
this, all functions are assumed to be smooth  in~$x$. We also
do not concern ourselves with sufficient  optimality conditions
(see, in this connection, for example,
\cite{car1,sagara,ssbook,tan_rugh}). Papers \cite{nnn}, \cite[\S
13]{kr_as} actually describe sufficient conditions of optimality for
same shadow price under sufficiently strong growth conditions.

 The structure of the paper is as follows: We begin with formulating
  the general control problem and stating general notation and  main
  assumptions 
  (Section \ref{sec:1}). Then, we formulate certain useful propositions from topology
   and stability theory
   (Section \ref{sec:11}) which are later used mostly in proofs;
     these propositions  are proved in Appendix.
    After that we discuss the relations of the Maximum Principle and
   introduce
     the notion of  $\tau$-vanishing Lagrange multipliers.
     Then we show that its existence is the necessary   optimality
     condition 
   ({Theorem} \ref{2190}). Connection between $\tau$-vanishing Lagrange multiplier and degenerate problems is investigated in Subsection \ref{ssect.zero};
for information on the connection with the condition $\psi^0(t)\to 0$
       refer to Subsection \ref{sec:23}.
      The problems with monotonic right-hand side  are investigated in
      Subsection \ref{sec:63}.
        Section \ref{sec:6} is mainly aimed at obtaining the most diverse
        sets
         of conditions  under which  a $\tau$-vanishing
         shadow price can be  explicitly expressed via a Cauchy-type formula.
Here we also discuss connections with the results of
\cite{kr_as,kab,av,norv}.

The last section is completely devoted to analysis of examples. We
show how the choice of a sequence of $\tau$ from a number of
uniformly weakly optimal solutions selects what is needed most
with the help of $\tau$-vanishing shadow price (Example
\ref{sternstern}).
 Example \ref{linlin} demonstrates that finding the $\tau$-vanishing Lagrange multiplier allows to solve abnormal problems in almost the same way as normal problems are solved. {Example} \ref{arnarn} shows how hard it is
to determine a $\tau$-vanishing Lagrange multiplier in cyclic
problems even if we know that the optimal control is unique. In
{Example}~\ref{avav}, the search for an optimal solution is reduced
to a boundary value problem.


  A part of results of this paper was announced in \cite{my},\cite{my1}.

\section{Preliminaries}
\label{sec:1} \ee
 We consider the time interval
  $\mm{T}\rav
  \mm{R}_{\geq 0}.$
The phase space of the
 control system is the
finite-dimensional metric space $\mm{X}\rav\mm{R}^m$; denote the
  unit ball in~$\mm{X}$  by~$\mm{D}$. Denote by $\mm{L}$
 the linear space of all real $m\times
m$ matrices; equip~$\mm{L}$ with the operator norm. The symbol~$E$
(which may be equipped with some indices) denotes various auxiliary
finite-dimensional Euclidean spaces.

For a subset~$A$ of a topological space,  denote by $cl\, A$ the
closure of $A$, and  by $int\, A$  the interior of $A$.
%
%

Slightly simplifying the notation when passing from the sequence
$\tau\rav(\tau_n)_{n\in\mm{N}}$ to its subsequence $\tau'$, we will
plainly write ``subsequence $\tau'\subset\tau$''

 Let $C(T,E)$ and $C_{loc}(T,E)$ be topological
spaces of all continuous functions of $T$ to $E$. Let us equip the
first one with extended  norm  $||\cdot||_C$ of  uniform
convergence. The second one is equipped the compact-open topology.

Here and below, for each integrable function~$a$ of time, the
integral $\int_0^\infty a(t)dt$ is the limit $\int_{0}^{T}a(t)dt$ as
$T\to\infty$. An improper integral, for example, over $[T,\infty\ra$,
is interpreted in the same sense.

Let us also consider a finite-dimensional Euclidean space~$\mm{U}$
and 
 map $U$ of $T$ to set of all subset of $\mm{U}$.   The
set~$\fr{U}$ of admissible controls is understood as the set of all
Borel measurable locally bounded selectors of the multi-valued
map~$U$. The topology on~$\fr{U}$ is defined through the inclusion
  $\fr{U}\subset \ct{L}^1_{loc}(\mm{T},\mm{U}).$

A function $a:\mm{T}\times E_1\times \mm{U}\to E_2$ is said to
\begin{description}
\item[1)] satisfy the Carath\'{e}odory conditions if
 a) the function $a(\cdot,x,u):\mm{T}\to E_2$ is Borel measurable for all
 $(x,u)\in E_1\times\mm{U},$
 b) the function $a(t,\cdot,\cdot):E_1\times \mm{U}\to E_2$ is continuous for a.a. $t\in\mm{T}.$
\item[2)]  be locally Lipshitz continuous
  if for each compact subset $K$ of $E_1\times\mm{U}$ there exists
  a function $L_K^a\in \ct{L}^1_{loc}(\mm{T},\mm{T})$
  satisfying
  $||a(t,x,u)-a(t,y,u) ||_{E_2}\leq L_K^a(t)||x-y||_{E_1}$ 
   for all
  $(x,u),(y,u)\in K$, $t\in\mm{T}$.
\item[3)]  be 
 integrally bounded (on each compact subset)
 if for each compact subset $K$ of $E_1\times\mm{U}$
  there exists a function
 $M_K^a\in \ct{L}^1_{loc}(\mm{T},\mm{T})$
 satisfying
  $||a(t,x,u)||_{E''}\leq M_K^a(t)$
  for all
  $(x,u)\in K$, $t\in\mm{T}$.
\end{description}

\medskip

 We assume the following conditions hold:

 {\it \bf Condition}~${\bf{(u)}}:$
    $U$ is a compact-valued map, and its
     graph is Borel set.

 {\it \bf Condition}~${\bf{(fg)}}:$   Locally Lipshitz continuous on $x$
 Carath\'eodory
 functions $f:\mm{T}\times\mm{X}\times\mm{U}\to\mm{X},$\
           $g:\mm{T}\times\mm{X}\times\mm{U}\to\mm{R},$\
           $\frac{\partial f}{\partial x}:\mm{T}\times\mm{X}\times\mm{U}\to\mm{L},$\
           $\frac{\partial g}{\partial x}:\mm{T}\times\mm{X}\times\mm{U}\to\mm{X}$
are integrally bounded
 (on each compact subset); in addition, $f$ satisfies the
 sublinear growth condition (see, for example, \cite[1.4.4]{tovst1}).


Let us
consider the control system
\begin{subequations}
   \begin{equation}
   \label{sys}
    \dot{x}=f(t,x,u),\  x(0)=x_{**},\qquad t\in\mm{T},\ x\in\mm{X},\ u(t)\in
    U(t),
   \end{equation}
 where $x_{**}\in\mm{X}$ is an initial value.
Now we can assign the solution of \rref{sys} to each~$u\in\fr{U}$.
The solution is unique and it can be extended to the whole $\mm{T}$.
Let
us denote it by~$x^u$. 
The map $u\mapsto x^u$ of $\fr{U}$ to $C_{loc}(\mm{T},\mm{X})$ is
continuous \cite{tovst1}.

In what follows, we study the problem of maximizing the objective
integral functional
\begin{equation}
   J^u(T)\rai\max;\qquad
   \label{opt}
   J^u(T)\rav\int_0^T g\big(t,x^u(t),u(t)\big)dt.
   \end{equation}
\end{subequations}
If there is no  limit in~\rref{opt}, the optimality may be defined in diverse ways
 (for details, see~\cite{car1,car,stern,tan_rugh});
generally, we will use the following definition:
\begin{definition}
We say that a control $u^0\in\fr{U}$ is {\it weakly uniformly
overtaking optimal} (see \cite{carm}) if
 $$ \limsup_{t\to\infty} \sup_{u\in \fr{U}} \big(J^u(t) - J^{u^0}(t)\big)
 =0.$$

For every sequence $\tau\rav(\tau_n)_{n\in\mm{N}}\uparrow\infty$ of
times, we say that a control $u^0\in\fr{U}$ is {\it $\tau$-optimal}
if
  $$\limsup_{n\to\infty} \sup_{u\in \fr{U}} \big(J^u(\tau_n) - J^{u^0}(\tau_n)\big) =0.$$
\end{definition}

We also assume:

  {\it \bf Condition~}${{(\tau)}}:$ there exists a weakly uniformly overtaking optimal control
   $u^0\in\fr{U}$ for problem \rref{sys}�-\rref{opt}.

 By this condition there exist an unbounded sequence
 $\tau=(\tau_n)_{n\in\mm{N}}\in \mm{T}^\mm{N}$ and some
 sequence $(\gamma_n)_{n\in\mm{N}}\in \mm{T}^\mm{N}$, converging to
zero,
 such that
 \begin{equation}
 \label{555}
{J}^{u^0}(\tau_n)\geq {J}^u(\tau_n)-\gamma^2_n\qquad\forall
u\in
 {\fr{U}},n\in\mm{N}.
   \end{equation}
   Then the control $u^0$ is $\tau$-optimal.
 Fix a sequence
$\tau$.  Also denote by~$x^0$ the trajectory that corresponds
to~$u^0$ .

Thus, any weakly uniformly  overtaking
 optimal control is {$\tau$-optimal} for some sequence
 $\tau\uparrow\infty$. Similarly,
 any
 uniformly overtaking  \cite{carm,2}       optimal
  control is  {$\tau$-optimal} for every sequence $\tau\uparrow\infty$.
 Since the definition of $\tau$-optimality refines these definitions, it is especially convenient if such sequence $\tau$ is given initially.

\section{{Auxiliary results}}
\label{sec:11}

\subsection{The set $\td{\fr{U}}$ of generalized
controls}  For each $u\in\mm{U}$, the symbol $\td\delta(u)$
denotes the probability measure concentrated at the point $u$.
 Denote by $\td{\fr{U}}_n$ the family of all weakly measurable
mappings~$\eta$ of $[0,n]$ to the set of Radon probability measures
over~$ \mm{U}$ such that $\eta(U(t))=1$ for a.a. $t\in[0,n]$. Let us
equip this set with the topology of *-weak convergence. Then, the
obtained topological space is a compact~\cite[IV.3.11]{va}, and the
set $\fr{U}_n\rav\{u|_{[0,n]}\,|\,u\in\fr{U}\}$ is everywhere densely
included in $\td{\fr{U}}_n$ \cite[IV.3.10]{va} by the map
$u\to\td{\delta}\circ u$.  We also keep the notation
$\td{u}^0\rav\td{\delta}\circ u^0.$

Now, let us introduce the set of all maps~$\eta$ of~$\mm{T}$ into the
set of Radon probability measures over~$ \mm{U}$ such that
$\eta|_{[0,n]}\in \td{\fr{U}}_n$ for every $n\in\mm{N}$, and let us
denote it by $\td{\fr{U}}$. For every~$n\in\mm{N}$, let the
projections $\td\pi_n:\td{\fr{U}}\to \td{\fr{U}}_n$  be given by
$\td{\pi}_n(\eta)\rav\eta|_{[0,n]}$ for all $\eta\in\td{\fr{U}}$. Let
us equip~$\td{\fr{U}}$ with the weakest topology such that all
projections are continuous. The set~$\td{\fr{U}}$ is called the set
of generalized controls.

\medskip

Let us assume that for the Euclidean space~$E$,  the function
 $a:\mm{T}\times{E}\times \mm{U}\to E$ is a locally Lipshitz continuous
  integrally bounded Carath\'{e}odory function that satisfies
   the extendability condition on~$\mm{T}$
  (for example, if the sublinear growth condition holds; see \cite[1.4.3]{tovst1}).

Let us fix a set  $\Xi\subset E$ of initial values and the system for
$u\in\fr{U}$:
\begin{subequations}
\begin{equation}
   \label{a}
   \dot{y}=a(t,y(t),u(t)),\  y(0)=\xi\in \Xi,\qquad
   t\in \mm{T}, u\in\fr{U}.
\end{equation}
It can also be generalized for $\eta\in\td{\fr{U}}$:
\begin{equation}
   \label{1650}
    \dot{y}=\int_{U(t)}a(t,y(t),u)d\eta(t),\ y(0)\in \Xi,\qquad
   t\in \mm{T}, \eta\in\td{\fr{U}}.
\end{equation}
 Each its local solution can be extended onto the
whole~$\mm{T}$. For every $\eta\in\td{\fr{U}}$, let us denote the
family of all solutions $y\in C_{loc}(\mm{T},E)$ of
system~\rref{1650} by $\td{\fr{A}}[\eta]$. Such transition from a
system defined for $u\in\fr{U}$ (like \rref{a}) to a generalized
system, which is defined for $\eta\in\td{\fr{U}}$ (like \rref{1650}),
will be done sufficiently often; to avoid writing the generalized
relation, we will write the initial one with the sign ``$\td{\ }$.''
For example, we will write $\td{\mbox{\rref{a}}}$ instead
of~$\mbox{\rref{1650}}$. In particular, for a solution $x^\eta\in
C_{loc}(\mm{T},\mm{X})$ of the Cauchy
problem~$\td{\mbox{\rref{sys}}}$, the function
  $T\mapsto
  \td{J}^\eta(T)$ could by introduced, for every $\eta\in \td{\fr{U}}$,
  by the rule $\td{\mbox{\rref{opt}}}$.

\begin{proposition}
 \label{spectr}
 Assume ${\bf{(u)}}$. Then,\
\begin{description}
\item[1)] the space~$\td{\fr{U}}$ is a compact, and $\td\delta({\fr{U}})$ is everywhere dense in it;
\item[2)]
   the map $\td{\fr{A}}:\td{\fr{U}} \to C_{loc}(\mm{T},E)$ is continuous and
     $\td{\fr{A}}[\td{\delta}\circ\fr{U}]$ is everywhere dense in a compact
      $\td{\fr{A}}[\td{\fr{U}}]\subset C_{loc}(\mm{T},E)$ for any compact $\Xi\subset E$;
\item[3)] If ${\bf{(fg)}}$ holds, then the map
  $\eta\mapsto x^\eta$ of $\td{\fr{U}}$ to $C_{loc}(\mm{T},\mm{X})$ and the map
  $\eta\mapsto \td{J}^\eta$
  of $\td{\fr{U}}$ to $C_{loc}(\mm{T},\mm{R})$
 are continuous.
\end{description}
\end{proposition}
 Since the proof of this proposition only plays an auxiliary role, it was repositioned to Appendix. Let us also note that embedding of the initial space $\fr{U}$ of admissible controls into a space with a more convenient topology is a well-known trick; see, for example,
 \cite{ga,va}, and  \cite{carm,dha,njota},
 \cite[Sect.~8]{kr_as} for infinite horizon problems.
 A weak  compactness  was used, for example, in   \cite{bald,car,dm,my2}.


\subsection{{Stability and thin tubes of solutions}} \label{sec:33}

 Let $w:\mm{T}\times \mm{U}\to\mm{T}$
  be an integrally bounded (on each compact subset) Carath\'{e}odory map. For
all $\tau\in\mm{T}$ and $\eta\in\td{\fr{U}}$, let us introduce
  $$\displaystyle \fr{L}_w[\eta](\tau)\rav\int_0^\tau\int_{U(t)} w(t,u)\,d\eta(t)\,dt.$$
Let us assume that $\fr{L}_w[\td{u}^0]\equiv 0$, and for every
$\eta\in \td{\fr{U}}$ from $\fr{L}_w[\eta](\tau)=0$ for all
$\tau\in\mm{T}$ it follows that $\eta$ equals $\td{u}^0$ a.e. on
$[0,\tau]$. The set of such~$w$ is denoted by $(Null)(u^0)$.

 For every position
$(\vartheta^*,y^*)\in\mm{T}\times{E}$, there exists a unique
solution~$y\in C(\mm{T},E)$ of the equation
\begin{equation}
   \label{1667}
    \dot{y}=a(t,y(t),u^0(t)),\quad y(\vartheta^*)=y^*.
\end{equation}
The solution
continuously depends on $(\vartheta^*,y^*)$. Let us denote its
initial position~$y(0)$ by $\beg(\vartheta^*,y^*)$.

\medskip

\begin{proposition}
\label{dop} 
     Let $U$ be a compact-valued
      map,
     and its
     graph is Borel set.
     Let $\Xi$ be a compact subset of
     $E.$

Then, there exists $w^0\in(Null)(u^0)$ 
such that for arbitrary $\eta\in\td{\fr{U}}$, $T\in\mm{T}$ for every
 solution
$y$ of \rref{1650} from $\beg(\vartheta,y(\vartheta))\in\Xi$ for all
$\vartheta\in
 [0,T]$ it follows that
 $$  ||\beg(\vartheta,y(\vartheta))-y(0)||_E\leq \fr{L}_{w^0}[\eta](\vartheta)\qquad\forall
 \vartheta\in
 [0,T].$$
 \end{proposition}
In the geometric sense, this proposition means that if a solution
$y|_{[0,T]}$ from the funnel~$\td{\fr{A}}[\eta]$ does not escape the
area $\td{\fr{A}}[u^0]$, then it also does not escape the tube of
   of  solutions of \rref{1667}, breadth of which (at $t=0$) does not surpass $\fr{L}_{w^0}[\eta](T)$. See the proof in Appendix.
\end{subequations}

\section{$\tau$-vanishing Lagrange multiplier as a necessary condition}
\label{sec:2}
\ee
\subsection{The core relations of the  Maximum Principle} \label{sec:21}

In what follows, we consider the shadow price~$\psi$ a covector (a
row vector); however, we will still write $x\in\mm{X},\psi\in\mm{X}$
and will not distinguish between the space $\mm{X}$ and its conjugate
space in the sense of sets.

  Let the Hamilton--Pontryagin function
  $\ct{H}:\mm{X}\times \mm{T}\times\mm{U}\times\mm{T}\times\mm{X}\to\mm{R}$
 be given by
 $$  \ct{H}(x,t,u,\lambda,\psi)\rav\psi f\big(t,x,u\big)+\lambda
 g\big(t,x,u\big).$$
Let us introduce the relations and boundary condition:
 \begin{subequations}
\begin{eqnarray}
       \dot{x}(t)&=& f\big(t,x(t),u(t)\big);\label{sys_x}
        \\
   \label{sys_psi}
       \dot{\psi}(t)&=&-\frac{\partial
       \ct{H}}{\partial x}\big(x(t),t,u(t),\lambda,\psi(t)\big);\\
   \label{maxH}
\sup_{p\in U(t)}\ct{H}\big(x(t),t,p,\lambda,\psi(t)\big)&=&
        \ct{H}\big(x(t),t,u(t),\lambda,\psi(t)\big);
   \end{eqnarray}
          \end{subequations}
  \begin{subequations}
   \begin{eqnarray}
      \label{dob}
     x(0)=x_{**},\quad       ||\psi(0)||_\mm{X}+\lambda&=&1.
   \end{eqnarray}
 It is easily seen that, for each $u\in\fr{U}$ for each initial
  condition, system~\rref{sys_x}--\rref{sys_psi} has a local solution,
   and each solution of these relations can be extended to the
    whole~$\mm{T}$. Let us denote by~$\fr{Y}$ the family of all solutions
   $(x,u,\lambda,\psi)\in
     C_{loc}(\mm{T},\mm{X})\times\fr{U}\times[0,1]\times C_{loc}(\mm{T},\mm{X})$
 of system~\rref{sys_x}--\rref{sys_psi},\rref{dob} on~$\mm{T}$.
 Let us denote by~$\fr{Z}$ the set of solutions from~$\fr{Y}$
 such that~\rref{maxH} also holds a.e. on~$\mm{T}$.

 Let us embed the sets ${\fr{Y}}$ and ${\fr{Z}}$ into $C_{loc}(\mm{T},\mm{X})\times\td{\fr{U}}\times[0,1]\times
 C_{loc}(\mm{T},\mm{X})$ by the mapping $(Id,\td{\delta},Id,Id)$; denote
  closures of their images by $\td{\fr{Y}}$ and $\td{\fr{Z}}$, respectively; then, $\td{\fr{Y}}$ and $\td{\fr{Z}}$ are compacts.

By Proposition \ref{spectr}, for every $(x,\eta,\lambda,\psi)\in
\td{\fr{Y}}$, the following relations hold: \rref{dob},
 $\td{\mbox{\rref{sys_x}}}$--$\td{\mbox{\rref{sys_psi}}}$; for $(x,\eta,\lambda,\psi)\in\td{\fr{Y}}$, we also have             $\td{\mbox{\rref{maxH}}}$, i.e.,
$$        \sup_{p\in
        U(t)}\ct{H}\big(x(t),t,p,\lambda,\psi(t)\big)=
        \int_{U(t)}\ct{H}\big(x(t),t,u,\lambda,\psi(t)\big)\,d\eta(t).
        \eqno{\tref{maxH}} 
        $$
Moreover, Proposition \ref{spectr} implies that all solutions of
these equations depend on both controls $u\in\fr{\td{U}}$ and initial
conditions continuously on any compact.

A nontrivial Lagrange multiplier  $(\lambda,\psi)\in
  [0,1]\times C_{loc}(\mm{T},\mm{X})$
   is called {\it{a  Lagrange multiplier associated with $(x^0,u^0)$}}
  if $(x^0,u^0,\lambda,\psi)$ is a solution of the core  Maximum
  Principle, i.e. the system
  \rref{sys_x}--\rref{maxH}.
It is convenient to denote by~$\Lambda$ the family of all Lagrange multipliers $(\lambda,\psi)\in
  \{0,1\}\times C_{loc}(\mm{T},\mm{X})$  associated with $(x^0,u^0)$ 
  such that
   \begin{eqnarray}
       \lambda=1\ or &\ & (\lambda=0\ and\ ||\psi(0)||_\mm{X}=1).
      \label{dobb}
   \end{eqnarray}
  \end{subequations}

\medskip

For each $\xi\in\mm{X}$, let us also define solutions $x_\xi\in C(\mm{T},\mm{X}),
 A_\xi\in C(\mm{T},\mm{L})$ of the following equations:
\begin{subequations}
  \begin{eqnarray}
   \dot{x}_\xi(t)=&\  &f(t,x_\xi(t),u^0(t))\quad x_\xi(0)=x_{**}+\xi, \\
   \label{sys_la}
         \dot{A}_\xi(t)=&\ &\frac{\partial
f}{\partial x}(t,x_\xi(t),u^0(t))\, A_\xi(t)\quad A_\xi(0)=1_\mm{L}.
   \end{eqnarray}
For every $T\in\mm{T}$, consider the covector
$$   I_{\xi}(T)\rav\int_0^T
   \frac{\partial g}{\partial x}(t,x_{\xi}(t),u^0(t))\, A_\xi(t)
  \,dt.$$

Similarly, for each $u\in{\fr{U}}$, let us  introduce a matrix
function $A^u$ and a covector function $I^u$ by the relations
\begin{eqnarray}
   \label{Aeta}
 \dot{A}^u(t)=\frac{\partial
f}{\partial x}(t,x^u(t),u(t))\,A^u(t),\qquad A^u(0)=1_\mm{L},\\
   I^{u}(T)\rav\int_0^T
   \frac{\partial g}{\partial x}(t,x^u(t),u(t))\, A^u(t)
  \,dt\quad \forall T\in\mm{T}.\nonumber
\end{eqnarray}
In addition, we call $x^\eta,A^\eta,\psi^\eta,I^\eta$  the solutions
of the corresponding $\td{\ }$-equations, or, equivalently, the limits, uniform on
compacts, of $x^u,A^u,\psi^u,I^u$ as $\td{\delta}(u)\to\eta$ in the
$*$-weak topology of $\td{\fr{U}}.$

 Expressing the solution of linear
equation \rref{sys_psi} through \rref{Aeta} (or \rref{sys_la}),
then any shadow price $\psi$ has the form
  \begin{equation}
   \label{cauchy}
         \psi(T)=(\psi(0)-\lambda I(T))A^{-1}(T)\quad \forall T\in\mm{T};
   \end{equation}
\end{subequations}
and  we can reformulate the result of \cite{Halkin} in the following way:
\begin{theorem}
\label{1}
    Under conditions $\bf{(u),(fg)}$,
for any $\tau$-optimal pair     $(x^0,u^0)\in C_{loc}(\mm{T},\mm{X})\times\fr{U}$ of problem \rref{sys}--\rref{opt}, for some  $\lambda^0\in [0,1],\psi^0\in C(\mm{T},\mm{X})$, the core relations of the Maximum Principle \rref{sys_x}--\rref{dob} hold for $(x^0,u^0,\lambda^0,\psi^0)$, i.e., $(x^0,u^0,\lambda^0,\psi^0)\in\fr{Z}.$

Moreover, up to a positive factor, for some
$I_*\in\mm{X},\iota_*\in\mm{X}$, one of the two following relations
also holds:
\begin{subequations}
  \begin{eqnarray}
   \label{h1}
     \lambda^0=1,&\  &\psi^0(T)=(I_*-I_0(T))A^{-1}_0(T)\qquad \forall T\in\mm{T};\\
   \label{h2}
    {\lambda}^0=0,&\ &\psi^0(T)=\ \ \ \ \ \ \ \ \ \ \ \, \iota_*A^{-1}_0(T)\qquad \forall T\in\mm{T}.
   \end{eqnarray}
\end{subequations}
\end{theorem}

 The core relations of the Maximum Principle are
incomplete, since \rref{sys_x}--\rref{dob}
do not contain a condition on the right endpoint, or, which is
actually equivalent, on $I_*$ or $\iota_*$. The remaining part of the
paper is mainly devoted
 to finding the additional relations at $I_*$ and $\iota_*$
  with the aid of $\tau$-vanishing Lagrange multiplier.


\subsection{Existence of $\tau$-vanishing multipliers}
            \label{sec:34}

System~\rref{sys_x}--\rref{sys_psi} can be rewritten for $u=u^0$ in
the form
\begin{subequations}
  \begin{eqnarray}
   \label{sys_psii}
       \dot{\psi}(t)=-\frac{\partial\ct{H}}{\partial
       x}(x(t),t,u^0(t),\lambda,\psi(t)),\\
   \label{sys_xxl}
       \dot{x}(t)=f(t,x(t),u^0(t)),\\
   \label{sys_lam}
       \dot{\lambda}=0.
   \end{eqnarray}
\end{subequations}

\begin{definition}

  A nontrivial Lagrange multiplier  $(\lambda^0,\psi^0)$ associated with $(x^0,u^0)$
   is called {\it $\tau$-vanishing}
  if $(\psi^0,x^0,\lambda^0)$ is a pointwise limit of a sequence of
 solutions $(\psi_n,x_n,\lambda_n)_{n\in\mm{N}}$ of system
  \rref{sys_psii}--\rref{sys_lam} such that $\psi_n(\tau'_n)=0$ for every $n\in\mm{N}$,
   here $\tau'\subset\tau.$
    In this case, the shadow price  $\psi^0$ is called {\it $\tau$-vanishing} as well.

\end{definition}
 Geometrically, this property means that the tube of solutions of system
 \rref{sys_psii}--\rref{sys_lam}, however thin (at the initial time), intersects
 with the hyperplane $\psi=0_\mm{X}$ at arbitrarily far time $\tau_n$.

 We claim that the existence of $\tau$-vanishing
 multipliers is a necessary optimality condition.
The main work horse of this proof is the following
    asymptotic condition of optimality  structurally
     similar to \cite[Theorem 9.1]{kr_as},\cite[Theorem 3]{kab}.

\begin{proposition}
 \label{1014}
Under conditions $\bf{(u),(fg)},(\tau)$, for each weight
$w\in(Null)(u^0),$
  there exist a sequence
  $(x^n,\eta^n,\lambda^n,\psi^n)_{n\in\mm{N}}\in \td{\fr{Y}}^\mm{N}$ and a
  subsequence $\tau'$ of
  $\tau$  such that
  \begin{description}
    \item[1)] for some $(x^0,\td{u}^0,\lambda^0,\psi^0)\in
   {\fr{Z}}$ it is $(x^n,\eta^n,\lambda^n,\psi^n)\to (x^0,\td{u}^0,\lambda^0,\psi^0)$
 in the topology of 
  $C_{loc}(\mm{T},\mm{X})\times\td{\fr{U}}\times[0,1]\times C_{loc}(\mm{T},\mm{X})$;
    \item[2)] $||\fr{L}_w(\eta^n)||_C\to 0$;
    \item[3)] $\td{J}^{\eta^n}(\tau'_n)-{J}^{u^0}(\tau'_n)\to 0+$; $\psi^n(\tau'_n)=0$ for all
$n\in\mm{N}.$
  \end{description}
\end{proposition}
The proof of this proposition  was repositioned to Appendix.

Note that from $\psi^n(0)=
   -\psi^n(\tau'_n)A^{\eta^n}(\tau'_n)+\psi^n(0)A^{\eta^n}(0)\ravref{cauchy}
   \lambda^n I^{\eta^n}(\tau'_n),$
   we have $\lambda^n I^{\eta^n}(\tau'_n)\to \psi^0(0).$


  Let $E=\mm{X}\times\mm{X}\times\mm{T}$,
$\Xi\rav 2\mm{D}\times(x_{**}+2\mm{D})\times  [0,1]$. To system
\rref{sys_psi},\rref{sys_x},\rref{sys_lam}, let us assign the
weight~$w$
 by means of Propositions~\ref{dop}.
Substituting this weight into Proposition~\ref{1014}, we obtain
\begin{remark}
 \label{1504}
Under conditions $\bf{(u),(fg)},(\tau)$
 there exist a sequence
  $(x^n,\eta^n,\lambda^n,\psi^n)_{n\in\mm{N}}\in \td{\fr{Y}}^\mm{N}$  and a
  subsequence $\tau'$ of
  $\tau$  such that
  \begin{description}
    \item[1)] for some $(x^0,\td{u}^0,\lambda^0,\psi^0)\in
   {\fr{Z}}$, it is $(x^n,\eta^n,\lambda^n,\psi^n)\to (x^0,\td{u}^0,\lambda^0,\psi^0)$
in the topology of
  $C_{loc}(\mm{T},\mm{X})\times\td{\fr{U}}\times[0,1]\times C_{loc}(\mm{T},\mm{X})$;
    \item[2)] the graphs of functions
   $(\psi^n,x^n,\lambda^n)$ are contained within the thinning funnels of solutions of
    system~\rref{sys_psii}--\rref{sys_lam}; i.e., for some sequence
    $(\delta_n)_{n\in\mm{N}}\in\mm{T}^\mm{N},$ $\delta_n\downarrow 0$, we have
   $$\varkappa(t,(\psi^n,x^n,\lambda^n)(t))\in (\psi^0(0),x_{**},\lambda^0)+
   \delta_n\mm{D}\times\delta_n\mm{D}\times[-\delta_n,\delta_n]
    \quad \forall t\in\mm{T},n\in\mm{N};$$
    \item[3)] $\td{J}^{\eta^n}(\tau'_n)-{J}^{u^0}(\tau'_n)\to 0+$;
    \item[4)] $\lambda^n I^{\eta^n}(\tau'_n)\to \psi^0(0)$; $\psi^n(\tau'_n)=0$ for all
$n\in\mm{N}.$
  \end{description}
\end{remark}
Note that $(\lambda^0,\psi^0)$ is nontrivial because it
  satisfies
  boundary condition \rref{dob} as well as the multipliers $(\lambda^n,\psi^n)$.
 For every $n\in\mm{N}$, consider a solution $(\psi_n,x_n,\lambda^n)$ of \rref{sys_psii}--\rref{sys_lam} with the initial conditions
 $(\psi_n(0),x_n(0),\lambda_n)\rav\varkappa(\tau'_n,(\psi^n(\tau'_n),x^n(\tau'_n),\lambda^n)).$
 Then $\psi_n(\tau'_n)=0_{\mm{X}}.$
 Since $(\psi_n(0),x_n(0),\lambda_n)\in(\psi^0(0),x^0(0),\lambda^0)+
   \delta_n\mm{D}\times\delta_n\mm{D}\times[-\delta_n,\delta_n],$
   i.e., $(\psi_n(0),x_n(0),\lambda_n)\to(\psi^0(0),x^0(0),\lambda^0),$
and because of the continuous dependency of solutions of \rref{sys_psii}--\rref{sys_lam},
    $(\lambda^0,\psi^0)$ is a $\tau$-vanishing Lagrange multiplier.
\begin{subequations}
 \begin{theorem}
 \label{2190}
 Assume that conditions $\bf{(u),(fg)},(\tau)$ hold.

Then, there exists a $\tau$-vanishing Lagrange multiplier
$(\lambda^0,\psi^0)\in\Lambda$, for example, constructed with a limit
of sequences from
   Remark~\ref{1504}.

Moreover, for every  $\tau$-vanishing Lagrange multiplier $(\lambda^0,\psi^0)\in\Lambda$,
  there exist
a subsequence $\tau'$ of
  $\tau$,
  a converging to~$0_\mm{X}$
 sequence $(\xi^n)_{n\in\mm{N}}\in \mm{X}^\mm{N}$,
 a converging to~$\lambda^0$ sequence $(\lambda^n)_{n\in\mm{N}}\in \la 0,1]^\mm{N}$
such that
\begin{eqnarray}
  \label{deimos_1}
   \psi^0(0)&=&\lim_{n\to \infty} \lambda^n I_{\xi^n}(\tau'_n);\\
  \psi^0(T)&=&
  \lim_{n\to \infty}\!
  \lambda^n(I_{\xi^n}(\tau'_n)-I_{\xi^n}(T))A_{\xi^n}^{-1}(T)
  \label{deimos_2}\\
  &=&\lim_{n\to
  \infty}\!
  \lambda^n\!\int_T^{\tau'_n}\!\frac{\partial g}{\partial x}(t,x_{\xi^n}(t),u^0(t))A_{\xi^n}(t)\,dt\,A_{0}^{-1}(T).
  \label{deimos_3}
\end{eqnarray}
and all the limits are uniform on every compact.

If, in addition to that, $\lambda^0>0$, then we can assume $\lambda_n=\lambda^0=1.$
\end{theorem}
\end{subequations}
\doc
The existence of a $\tau$-vanishing Lagrange multiplier
$(\lambda^0,\psi^0)$ is shown above. By multiplying this nontrivial
$(\lambda^0,\psi^0)$ by a certain scalar, we can always provide
condition~\rref{dob}; thus, $(\lambda^0,\psi^0)\in\Lambda$.

 Let $(\lambda^0,\psi^0)$ be  a $\tau$-vanishing Lagrange multiplier.
 The sequences $\tau',(\lambda_n)_{n\in\mm{N}},(\psi_n)_{n\in\mm{N}}$
exist by the definition of a $\tau$-vanishing Lagrange multiplier if
we define $\xi^n\rav x_n(0)-x^0(0)$ for every $n\in\mm{N}$. Since for
all $n\in\mm{N}$ $\psi_n(\tau'_n)=0_\mm{X}$, the Cauchy formula
\rref{cauchy}
 implies
  $\psi_n(T)=\lambda^n(I_{\xi^n}(\tau'_n)-I_{\xi^n}(T))A_{\xi^n}^{-1}(T)$      for every $T\in\mm{T},$
we have $\psi_n(0)=\lambda^n (I_{\xi^n}(\tau'_n)-I_{\xi^n}(T))=\lambda^n I_{\xi^n}(\tau'_n).$
Now, uniformity of the limit $\psi^0$ of $\psi_n$ yields \rref{deimos_1}.
Substituting this into \rref{cauchy} we obtain \rref{deimos_3}
for every $T\in\mm{T}$. What remains follows from the theorem of
continuous dependence of solutions on initial conditions, applied to
\rref{sys_psii}--\rref{sys_lam} and \rref{sys_la}. \bo

\subsection{On different topologies for the set of generalized controls} \label{sec:39}

 Consider a weight $w^0\in (Null)(u^0)$. Define $w^1$ by the rule
 $w^1(t,u)\rav w^0(t,u)+||u-u^0(t)||$ for every $(t,u)\in\mm{T}\times\mm{U}.$
 Then, for a subsequence $(u_n)_{n\in\mm{N}}\in\fr{U}^{\mm{N}}$,
 from $||\fr{L}_{w^1}[\td{\delta}\circ u_n]||_{C}\to 0$ it follows that
 $||u_n-u^0||_{\ct{L}^{1}(\mm{T},\mm{U})}\to 0$
  (certainly, this does not imply that
  $u^0\in\ct{L}^{1}(\mm{T},\mm{U})$). Similarly, for any
  $p\in\la 0,\infty\ra,\nu\in B_{loc}(\mm{T},\mm{R}_{>0})$,
replacing $||u-u^0(t)||$ with
  $\nu(t)||u-u^0(t)||^p$
guarantees the convergence of $u_n-u^0\to 0$ in the topology  of
  $\ct{L}^{p}_{\nu}(\mm{T},\mm{U}).$
For every interval $\fr{T}\subset\mm{T}$, this extended metric also
induces the extended distance
$\varrho\big(\eta,u^0;\td{\ct{L}}^{p}_{\nu}(\fr{T},\mm{U})\big)$ on
$\td{\fr{U}}$ by the rule
\begin{equation*}
  \varrho\big(\eta,u^0;\td{\ct{L}}^{p}_{\nu}(\fr{T},\mm{U})\big)
  \rav
  \Big(\int_{\fr{T}}\nu(t)\int_{U(t)}||u-u^0(t)||^p\,d\eta(t)\,dt
  \Big)^{1/p}\
   \forall \eta\in \td{\fr{U}}. 
\end{equation*}

Addition of the summand $\nu(t)R^p(t,u)$ (see~\rref{nado}) provides the uniform convergence
   $||\dot{y}(t)-a(t,y(t),u^0(t))||_{\ct{L}^{p}_{\nu}(\mm{T},\mm{X})}\to 0$ by all
  $\eta\in\td{\fr{U}},y\in\fr{A}[\eta]$ such that $y(0)\in \Xi$.

 Let us replace the weight $w$ from Proposition~\ref{1014} and
    {Remark} \ref{1504} by stronger ones if necessary. Then
    there exists a $\tau$-vanishing Lagrange multiplier
$(\lambda^0,\psi^0)$ as the limit of sequences from
   Remark~\ref{1504}.
Now we have
\begin{remark}
\label{sas}
  Assume that conditions $\bf{(u),(fg)},(\tau)$ hold.
Then there exists a $\tau$-vanishing multiplier $(\lambda^0,\psi^0)$
   associated with $(x^0,u^0)$ such that for this multiplier,
   the conclusion of {Remark} \ref{1504} holds and, moreover,
     the following convergences are guaranteed:
$\varrho\big(\eta^n,u^0;\td{\ct{L}}^{p}_{\nu}(\mm{T},\mm{U})\big)\to 0,$
  $||\dot{x}^n(t)-f(t,x^n(t),u^0(t))||_{\ct{L}^{p}_{\nu}(\mm{T},\mm{X})}\to 0_\mm{X}.$
\end{remark}
\medskip

The condition~${\bf (u)}$ implies that, a.a. $t\in\mm{T}$,
the controls are chosen from the compact~$U(t)$. Let us weaken this
assumption to the following:

{\it \bf Condition}~${\bf{(u_\sigma)}}:$
      $U$ is a closed-valued
     map, and its
     graph is Borel set.

We shall still assume the conditions ${\bf{(fg)}},(\tau)$ to hold.
A nondecreasing sequence $(U^{(r)})_{r\in\mm{N}}$ of locally bounded
 compact-valued maps  is given by
 $$U^{(r)}(t)\rav\{u\in U(t)\,|\,||u-u^0(t)||<r\}\quad \forall t\in\mm{N},r\in\mm{N}. $$
Let the set~$\fr{U}^{(r)}$ 
be the set of all Borel measurable selectors of the multi-valued
map~$U^{(r)}$.
 Then  for all $r\in\mm{N}$ $ u^0\in \fr{U}^{(r)}\subset \fr{U}^{(r+1)}$ and
  $\displaystyle U\equiv\cup_{r\in\mm{N}}U^{(r)}$ hold; now, we have
  $\displaystyle \fr{U}^{(\infty)}\rav\cup_{r\in\mm{N}}\fr{U}^{(r)}\equiv\fr{U}$.


Repeating the reasonings of~Sect~\ref{sec:11}, for every
$r\in\mm{N}\cup\{\infty\}$, we can construct sets $\td{\fr{U}}^{(r)}$
and images
  ${\fr{U}}^{(r)}_n\rav\pi_n(\fr{U}^{(r)}),
  \td{\fr{U}}^{(r)}_n\rav\td{\pi}_n(\td{\fr{U}}^{(r)}).$
  Denote by $\td{\fr{U}}$ the set   of all maps~$\eta$ from~$\mm{T}$ into
 the set of Radon probability measures over~$ \mm{U}$ such that
 $\eta|_{[0,n]}\in \td{\fr{U}}^{(\infty)}_n$  for every $n\in\mm{N}$. The topology
 of this set is the weakest topology in which $\fr{U}^{(r)}$ could
 be continuously embedded into $\td{\fr{U}}$.
 Note that under our definition, 
 $\td{u}^0\in \td{\delta}(\fr{U}^{(r)})\subset 
\td{\delta}(\fr{U})$ for all $r\in\mm{N}$.

To system \rref{sys_psi},\rref{sys_x},\rref{sys_lam}, let us assign
the weight~$w$
 by means of Propositions~\ref{dop}. Note that this weight is independent of $r.$
For the sequence $\tau$, for each
  $\td{\fr{U}}^{(r)}$, we have
  {Remark}
 \ref{1504}; in particular, there exist a time $t_r\in\tau,t_r>r,$ a
 $\tau$-vanishing Lagrange multiplier $(\lambda^r,\psi^r)$, and a solution
 $(x^r,\eta^r,\td{\psi}^r,\td{\lambda}^r)\in\td{\fr{Y}}$ with the properties
 \begin{subequations}
 \begin{eqnarray}\label{maxHr}
\sup_{p\in U^{(r)}(t)}\ct{H}\big(x(t),t,p,\lambda^r,\psi^r(t)\big)=
        \ct{H}\big(x(t),t,u^0(t),\lambda^r,\psi^r(t)\big) \forall
        \ a.a.\ t\in\mm{T}\\
         ||\fr{L}_w(\eta^r)||_C<1/r,\label{Hr1}
||\td{x}^r-x^r||_{C([0,r],\mm{X})}<1/r,\
   ||\td{\psi}^r-\psi^r||_{C([0,r],\mm{X})}<1/r,\\
   \big|\big|\varkappa(t_r,(\td{\psi}^r(t_r),\td{x}^r(t_r),\td{\lambda}^r))-
   (\psi^0,x_{**},\lambda^0)\big|\big|_{E}<1/r,
   \label{Hr2}\\
     0\leq\td{J}^{\eta^r}(t_r)-{J}^{u^0}(t_r)<1/r,\  \td{\psi}^r(t_r)=0_\mm{X}.
    \label{Hr3}
   \end{eqnarray}
\end{subequations}
Passing to the limit, we obtain $\eta^r\to\td{u^0}$ from
$||\fr{L}_w(\eta^r)||_C<1/r$. Passing to the subsequence
$\tau'\subset (t_r)_{r\in\mm{N}}\subset\tau$, we can provide the
monotonicity of $t_r$ and convergence of the sequence
   $(\lambda^r,\psi^r)_{r\in\mm{N}}\in\big(\la0,1]\times C_{loc}(\mm{T},\mm{X})\big)^\mm{N}$ to certain
   $(\lambda^0,\psi^0)$. Under these assumptions, we immediately see that
    $(\psi^0,x^0,\lambda^0)$ is the solution of
   {{\rref{sys_psii}}}{-}{{\rref{sys_lam}}} that satisfies
   \rref{dob}. Now the sequence
   $(\td{x}^r,\eta^r,\td{\psi}^r,\td{\lambda}^r)_{r\in\mm{N}}$ converges,
    by \rref{Hr1}, to $(x^0,\td{u}^0,\lambda^0,\psi^0).$
    Passing to the pointwise limit in \rref{maxHr}, we obtain for
   $(x^0,\td{u}^0,\lambda^0,\psi^0)$ the property \rref{maxH}.
   Thus we proved items 1) and 2) of {Remark} \ref{1504}. Since the limit of
  $(\td{x}^r,\eta^r,\td{\psi}^r,\td{\lambda}^r)_{r\in\mm{N}}$ and
  $({x}^0,\eta^r,{\psi}^r,{\lambda}^r)_{r\in\mm{N}}$ is the same, items 3) and 4) follow
  from \rref{Hr2} and \rref{Hr3} respectively.

Consider again the solutions $(\psi_n,x_n,\lambda^n)$ of
\rref{sys_psii}--\rref{sys_lam} for the initial conditions
 $(\psi_n(0),x_n(0),\lambda_n)\rav
 \varkappa(\tau'_n,(\td{\psi}^n(\tau'_n),\td{x}^n(\tau'_n),\td{\lambda}^n))$. Then
 $(\lambda^0,\psi^0)$ is a $\tau$-vanishing Lagrange multiplier and {Theorem}
 \ref{2190} holds under condition $\bf(u_\sigma)$. Thus,
\begin{corollary}
\label{sigma}
Condition  $\bf(u)$ in
  {Remark} \ref{1504},
  {Theorem}
 \ref{2190} could be replaced with
   $\bf(u_\sigma)$.
\end{corollary}

\begin{corollary}
 Assume conditions $\bf{(u_\sigma),(fg)}$ hold.
 Let a pair $(x^0,u^0)\in C_{loc}(\mm{T},\mm{X})\times\fr{U}$ be
  weakly uniformly overtaking  optimal
  for problem~\rref{sys}--\rref{opt}.

Then, for some unbounded sequence
 $\tau=(\tau_n)_{n\in\mm{N}}\in\mm{T}^\mm{N}$, there exists a $\tau$-vanishing Lagrange multiplier
$(\lambda^0,\psi^0)\in\Lambda$.
\end{corollary}
\begin{corollary}
 Assume conditions $\bf{(u_\sigma),(fg)}$ hold.
 Let a pair $(x^0,u^0)\in C_{loc}(\mm{T},\mm{X})\times\fr{U}$ be
  uniformly overtaking  optimal
  for problem~\rref{sys}--\rref{opt}.

Then, for each unbounded sequence
 $\tau=(\tau_n)_{n\in\mm{N}}\in\mm{T}^\mm{N}$, there exists a $\tau$-vanishing Lagrange multiplier
$(\lambda^0,\psi^0)\in\Lambda$.
\end{corollary}

\section{Properties of $\tau$-vanishing Lagrange multipliers}
\subsection{On  stable shadow prices}
\label{sec:23}

 Consider the boundary conditions
\begin{subequations}
\begin{eqnarray}
   \label{trans}
       \displaystyle \lim_{t\to\infty} \psi(t)=0,\\
   \label{partlim}
       \liminf_{n\to\infty}||\psi^0(\tau_n)||_{\mm{X}}=0.
   \end{eqnarray}
\end{subequations}
\begin{definition}
The component~$\psi^0$ of a solution~$y^0=(\psi^0,x^0,\lambda^0)$
 of system
    \rref{sys_psii}--\rref{sys_lam} is said to be {\it Lyapunov stable   in
domain~$\Xi$} if for each $\epsi>0$ there exists $\delta>0$ such that
for each solution~$y=(\psi,x,\lambda)$
 of system
    \rref{sys_psii}--\rref{sys_lam}
 from  $||y(0)-y^0(0)||_{E}<\delta, y(0)\in \Xi$ it follows
  that $||\psi^{0}(s)-\psi(s)||_\mm{X}<\epsi$ for all $s\in\mm{T}$.
\end{definition}
\begin{corollary}
\label{s3}
 Assume that conditions $\bf{(u_\sigma),(fg)},(\tau)$ hold.
   Let for some  solution $(\psi,x^0,\lambda)$
 of system
    \rref{sys_psii}--\rref{sys_lam}
the component $\psi$ be  Lyapunov stable in the domain
$\mm{X}\times\mm{X}\times[0,1]$.

 Then all
 $\tau$-vanishing multipliers $(\lambda^0,\psi^0)\in\Lambda$
satisfy the condition \rref{partlim}.
\end{corollary}
  \doc
Since equation
 \rref{sys_psii} is linear, the Lyapunov stability of $\psi$ for some solution $(\psi,x^0,\lambda)$
 of system
    \rref{sys_psii}--\rref{sys_lam} yields the Lyapunov stability of this component for all solutions of system
    \rref{sys_psii}--\rref{sys_lam}.

Consider  every
 $\tau$-vanishing multiplier
 $(\lambda^0,\psi^0)$ and the sequences $\tau',(\lambda_n)_{n\in\mm{N}},(x_n)_{n\in\mm{N}},(\psi_n)_{n\in\mm{N}}$
from its definition. Then, $y_n=(\psi_n(0),x_n(0),\lambda_n)\to
y^0=(\psi^0(0),x^0(0),\lambda^0),$ and by definition of Lyapunov
stability for some  $N\in\mm{N}$ for all $n\in\mm{N}, n>N$
   $||\psi^0(\tau'_n)||_\mm{X}=||\psi^n(\tau'_n)-\psi^0(\tau'_n)||_\mm{X}<\epsi$.
   Since $\epsi> 0$ was arbitrary, we have shown \rref{partlim}
for all $\tau$-vanishing multipliers.\bo


 Note that since~\rref{sys_psi} is linear, the partial stability of the variable~$\psi$
 implies its boundedness. Therefore, the proved proposition is useless if all shadow prices are unbounded.
Note that, as follows from \cite[Example~5.1]{stern},
 for a weakly uniformly overtaking optimal control $u^0$, a
 $(x^0,u^0,\lambda^0,\psi^0)\in\fr{Z}$ that satisfies \rref{partlim} may not satisfy stronger condition~\rref{trans}.

The stability condition can be replaced with a condition 
 which is stronger but much easier to check. 
 \begin{corollary}~\label{partlimz}
 Assume that conditions $\bf{(u_\sigma),(fg)},(\tau)$ hold.
 If
  the functions  $L_K^f,L_K^g$ are independent of the
compact~$K$, and
   these functions
   are
 summable
 on  $\mm{T}$ (see \cite[Hypotesis 3.1 (iv)]{sagara}), then
  any $\tau$-vanishing multiplier satisfies condition
 \rref{trans}.
\end{corollary}
%
\doc
 Let
$(\psi^0,\lambda^0)$ be a $\tau$-vanishing multiplier. Let $\xi_0\rav
(\psi^0(0),x^0(0),\lambda^0),$
$\Xi\rav\xi_0+\mm{D}\times\mm{D}\times[-1,1].$
 By
\cite[(3.3)]{sagara} there exists a summable function $\omega\in
\fr{L}^1(\mm{T},\mm{T})$
 such that $\dot{\psi}(t)\leq \omega(t)$ for a.a. $t\in\mm{T}$ for all
 solution $(\psi,x,\lambda)$
 of system
    \rref{sys_psii}--\rref{sys_lam}
    if $\xi\rav(\psi(0),x(0),\lambda)\in\Xi$.
     Now for each pair $(t_1,t_2)\in\mm{T},t_1\leq t_2,$
\begin{eqnarray*}
 ||\psi(t_1)-\psi^0(t_2)||_\mm{X}\leq||\psi-\psi^0||_{C([0,t_1],\mm{X})}+2
 \int_{t_1}^{\infty}\omega(t)\,dt
\end{eqnarray*}
 if $\xi\in\Xi$.
For each $\epsi>0$ there exists $T\in\mm{T}$ such that the second
summand does not exceed $\epsi/2$ if $t_1>T;$ now
 there exists $r\in\mm{T}$ such
that $||\psi-\psi^0||_{C([0,T],\mm{X})}$ does not exceed $\epsi/2$
 if $||\xi-\xi^0||_E<r.$
 Then, setting $t_1=t_2$, we obtain
 $||\psi-\psi^0||_C\leq\epsi$ if $||\xi_1-\xi_2||_E<r,$
 i.e.,
the component $\psi^0$ is Lyapunov stable. By Corollary~\ref{s3},
\rref{partlim} holds, and $||\psi^0(t_1)||<\epsi/2$  for some
$t_1\in\mm{T}, t_1>T$.

Then, setting $\xi=(\psi^0(0),x^0(0),\lambda^0)$, we obtain
 $||\psi^0(t_2)||_\mm{X}=||\psi^0(t_2)-\psi^0(t_1)||_\mm{X}+||\psi^0(t_1)||_\mm{X}<\epsi$
  if $t_2>t_1.$ Thus \rref{trans} holds.
\bo

The even more strong conditions used for proving transversality
condition~\rref{trans}  can be seen, for example,
 in~\cite[(A3)]{Ye} (the Lipshitz
constants  $L_K^g, L_K^f$ are required to decrease exponentially
with time).

\subsection{Degenerate $\tau$-vanishing Lagrange multipliers}
\label{ssect.zero}



\begin{remark}
 \label{2351111}
 Assume that conditions $\bf{(u_\sigma),(fg)},(\tau)$ hold.
  If  for some $\tau'\subset\tau$
\begin{equation}
 \label{1846}
\limsup_{n\to\infty, \xi\to 0_\mm{X}}||I_\xi(\tau'_n)||_\mm{X}
     <\infty
\end{equation}
then    the pair $(x^0,u^0)$ is normal,  there exists a
$\tau$-vanishing multiplier $(1,\psi^0)\in\Lambda$.

  Moreover, if  
$  \displaystyle \limsup_{n\to\infty, \xi\to
0_\mm{X}}||I_\xi(\tau_n)||_\mm{X}
     <\infty,$
then  every $\tau$-vanishing multiplier
$(\lambda^0,\psi^0)\in\Lambda$ satisfies
$\lambda^0=1.$

 On the other hand, if  
$  \displaystyle \lim_{ n\to\infty, \xi\to
0_\mm{X}}||I_\xi(\tau_n)||_\mm{X}
     =\infty,$
then  every $\tau$-vanishing  multiplier
$(\lambda^0,\psi^0)\in\Lambda$  satisfies
$\lambda^0=0.$
 \end{remark}
\doc By Theorem~\ref{2190} there exists a $\tau$-vanish multiplier
$(\lambda^0,\psi^0)\in\Lambda$ satisfying \rref{deimos_1}, but for
each such multiplier  we have
$\lambda^n||I_{\xi_n}(\tau'_n)||_{\mm{X}}=||\psi^n(0)||_{\mm{X}}
 \to ||\psi^0(0)||_{\mm{X}};$ then
 $\lambda^0=0$ iff $(||I_{\xi_n}(\tau'_n)||_{\mm{X}})_{n\in\mm{N}}\uparrow
 \infty.$ \bo

There are many conditions that provide nondegeneracy of the problem; in connection with this, note papers
 \cite{kr_as,kab,av,aucl,norv}.
The connection between the normality of the problem and finiteness
of~$I_0$ seems to be noted for the first time in
\cite[(3.24)]{kr_as}. Condition~\rref{1846} develops this approach,
actually demanding $I_\xi$ to be locally bounded. As we are going to
show below, many sufficient
 conditions of nondegeneracy for the optimal problem can be obtained from
 \rref{1846}.
However, there are other ways to prove the nondegeneracy. For example,
 \cite[Theorem~5.1]{kr_as} uses the smoothness of the objective value function, and
 \cite[Theorem~10.1]{kr_as} and \cite[Theorem~5]{kab} use the monotonicity of the functions $f$ and $g$ in $x$ and the stationarity condition.

\medskip

Note that although the examples of
 abnormal problems are well known (\cite{Halkin,kr_as,Pickenhain}),
additional relations of the Maximum Principle for such problems did
not receive much attention from researchers; the author only knows of
the dual problem construction in paper \cite{Pickenhain}. Let us
apply {Theorem}
 \ref{2190} to these problems.

\medskip

Consider a degenerate $\tau$-vanishing solution
$(x^0,u^0,0,\psi^0)\in\fr{Z}$. Then, from \rref{dob} we have
$\psi^0(0)=1$, and Theorem \ref{2190} yields
\begin{eqnarray}
   \label{1782}
  \psi^0(0)\ravref{deimos_1}\frac{\psi^0(0)}{||\psi^0(0)||_\mm{X}}=
  \lim_{n\to\infty}\frac{\lambda_n I_{x_n(0)}(\tau'_n)}{||\lambda_n I_{x_n(0)}(\tau'_n)||_\mm{X}}=
  \lim_{n\to\infty}\frac{I_{x_n(0)}(\tau'_n)}{||I_{x_n(0)}(\tau'_n)||_\mm{X}}
\end{eqnarray}
provided $x_n(0)\to x^0(0).$ Using Remark~\ref{2351111}, we finally
obtain

\begin{corollary}
\label{2353} Let $\bf{(u_\sigma)}$,$\bf{(fg)}$, $(\tau)$ hold. Let
  $$    \lim_{n\to\infty, \xi\to 0_\mm{X}}
  ||I_\xi(\tau_n)||_{\mm{X}}=\infty,\quad
        \lim_{n\to\infty, \xi\to 0_\mm{X}} \frac{I_\xi(\tau_n)}{||I_\xi(\tau_n)||_{\mm{X}}}=\iota_*.
  $$
 for some
vector $\iota_*\in\mm{X}$.

Then, there is unique $\tau$-vanishing multiplier
$(0,\psi^0)\in\Lambda$, and $\iota_*$ and $\psi^0$ are connected  by
\rref{h2}.
\end{corollary}

\subsection{Monotonic case}
 \label{sec:63}

Consider a nonempty convex closed cone $\fr{C}.$
 The cone orderings $\succcurlyeq,\succ$
   of $\mm{X}$  induced by~$\fr{C}$ are the relations
   defined as follows: for all $x,
   y\in\mm{X},\  $
   $$(x \succcurlyeq_\fr{C} y) \Leftrightarrow (x - y \in \fr{C}),\qquad
   (x \succ_\fr{C} y) \Leftrightarrow (x - y \in int\, \fr{C}).$$
The pre-orders  on $\mm{L}$
 defined as follows: for $B,C\in\mm{L},\  $
\begin{eqnarray*}
 (B\succcurlyeq_\fr{C} C) &\Leftrightarrow&
  ((B-C)x
  \in\fr{C}\quad\ \ \ \ \, \forall x\in\fr{C}),\\
 (B \succ_\fr{C} C) &\Leftrightarrow&
   ((B-C)x
   \in int\,\fr{C}\quad \forall x\in int\,\fr{C}).
\end{eqnarray*}
%
   Note that $1_\mm{L}\succcurlyeq_\fr{C} 0_\mm{L},$
   $1_\mm{L}\succ_\fr{C} 0_\mm{L}.$

   The  conjugate  cone  of $\fr{C}$
   is defined by  $\fr{C}^\bot\rav\{x\in\mm{X}\,|\,\forall y\in\fr{C}\ xy\geq 0\}$.

\begin{proposition}
\label{mon10}
 Assume that conditions $\bf{(u_\sigma),(fg)},(\tau)$ hold.
  Assume that there exists Carat\'{e}odory function $d:\mm{T}\times\mm{X}\to\mm{R}$
  such that for all $x\in\mm{X}$ and a.a. $t\in\mm{T}$
   the following relation holds:
 $$\frac{\partial g}{\partial x}(t,x,u^0(t))\succcurlyeq_{\fr{C}^\bot}
0_\mm{L},\quad
  \frac{\partial f}{\partial x}(t,x,u^0(t))\succcurlyeq_{\fr{C}}
  d(t,x)1_\mm{L}.$$
Then, there exists a $\tau${-vanishing} multiplier
 $(\lambda^0,\psi^0)\in\Lambda$,
  {and for any such  multiplier, we have} 
   $\psi^0\succcurlyeq_{\fr{C}^\bot} 0_\mm{X},$
  and $\psi^0(0)\in\fr{C}^\bot.$

Moreover, if $\lambda^0>0$ (for example, if \rref{1846} holds), then
for all $y\in\fr{C}$
\begin{equation}
 \label{2455}
   \limsup_{t\to\infty,\xi\to 0_\mm{X}} I_\xi(t)y
    \geq\psi^0(0)y\geq
   \lim_{t\to\infty} I_0(t)y \geq
   0,
\end{equation}
 and all limits in \rref{2455} are correctly defined.

 If, in addition, there exists a Lebesgue point
$t^*\in\mm{T}$ for the function
   $u^0$ such that
$$  \frac{\partial g}{\partial x}(t^*,x^0(t^*),u^0(t^*))\succ_{\fr{C}^\bot} 0_\mm{L},$$
then
 $
   \psi^0|_{[0,t^*]}\succ_{\fr{C}^\bot} 0_\mm{X};$
in particular, $\psi^0(0)\in int\, \fr{C}^\bot$. If such Lebesgue
point $t^*$ exists  on every infinite interval,
 then $\psi^0\succ_{\fr{C}^\bot} 0.$
\end{proposition}
%
%
%
\doc
Fix arbitrary $\xi\in\mm{X}$, $T>0$,
   $\vartheta>T$.
%
   Denote by $F_\xi(t)$ the matrix
   $\frac{\partial f}{\partial x}(t,x_\xi(t),u^0(t))$, and
    by $m_\xi$ the measurable function
   $t\mapsto -d(t,x_\xi(t))$;
   by condition, $F_\xi+m_\xi(t)1_\mm{L}
   \succcurlyeq_{\fr{C}} 0_\mm{L}.$ Now, let us consider a solution~$P(t)$ of the equation
      $$\dot{P}=\big(F_\xi(t)+m_\xi(t)1_\mm{L}\big) P,\quad P(T)=1_\mm{L},\quad t\geq T;$$
       then $P(t) \succcurlyeq_{\fr{C}} 1_\mm{L}$ for all $t\in\la
T,\vartheta].$ But the solution $P$ is the product of two nonnegative
solutions of the equations
      $\dot{Q}=F_\xi(t)Q,\ Q(T)=1_\mm{L}$, and
      $\dot{r}_\xi=m_\xi(t)r_\xi,\ r_\xi(T)=1.$
Thus, $P(\vartheta)=r_\xi(\vartheta)Q(\vartheta)=
      r_\xi(\vartheta)A_\xi(\vartheta)A^{-1}_\xi(T),$
and $P(\vartheta)\succcurlyeq_{\fr{C}} 1_\mm{L}$ implies
$A_\xi(\vartheta)A^{-1}_\xi(T)=Q(\vartheta)=P(\vartheta)/r_\xi(\vartheta)\succcurlyeq_{\fr{C}}
1_\mm{L}/r_\xi(\vartheta) $ for all $\vartheta >T$. In particular,
for all $y\in\fr{C}$, we have
 $A_\xi(\vartheta)A^{-1}_\xi(T)y \succcurlyeq_{\fr{C}} y/r_\xi(\vartheta) $,
whence
\begin{eqnarray}
   \ \!\!\!\!\!\!\!\!\!\frac{dI_\xi(t)}{dt}A_\xi^{-1}(T) y = \!
   \frac{\partial g}{\partial
    x}(t,x_\xi(t),u^0(t)) A_\xi(t)A_\xi^{-1}(T)y
    \!
    \geq\! \frac{\partial g}{\partial
    x}(t,x_\xi(t),u^0(t)) \frac{y}{r_\xi(t)}
    \!\geq\! 0    \label{mmm}
\end{eqnarray}
 for all $\xi\in\mm{X},y\in\fr{C},T\in\mm{T},t>T;$
thus, for $T=0$, we have $\frac{dI_\xi(t)}{dt} \in \fr{C}^\bot$,
hence  the functions
    $I_\xi y,I_\xi  A_\xi^{-1}(T)y$ are monotonic
 for all
 $\xi\in\mm{X},T\in\mm{T},y\in\fr{C}.$

By Theorem~\ref{2190}, there exists a $\tau$-vanishing multiplier
 $(\psi^0,\lambda^0)\in\Lambda$. Moreover,
 each such multiplier
 $(\psi^0,\lambda^0)\in\Lambda$
   satisfies
formula~\rref{deimos_3} for  certain sequences $\lambda^n$ and
$\xi_n$. However, the
 integrand of~\rref{deimos_3}
    lies in $\fr{C}^\bot$. Passing to the limit as $n\to\infty$,
     we obtain $\psi^0\succcurlyeq_{\fr{C}^\bot} 0_\mm{X}.$

Fix the basis of $span\, \fr{C}$ made of the vectors $y\in\fr{C}$;
now, for every such vector $y$, the functions
    $I_\xi y$ 
    are monotonic,
and
$$  \limsup_{t\to \infty,\xi\to 0_\mm{X}}{\lambda_0} I_{\xi}(t)y\geq
  \lim_{n\to \infty}  {\lambda_0}
I_{\xi_n}(\tau'_n)y
 \ravref{deimos_1}{\psi^0(0)}y;$$
we obtain the first estimate from \rref{2455}.


   Fix any $T\in\mm{T},y\in\fr{C}.$ 
  Now, monotonicity of~$I_\xi A_\xi^{-1}(T) y$ 
 yields
\begin{eqnarray}
\psi^0(T)y\!\ravref{deimos_2}\!\lim_{n\to \infty}\lambda^n
  \big(I_{\xi_n}(\tau'_n)-I_{\xi_n}(T)\big)
  A_{\xi_n}^{-1}(T)y\geq \lim_{n\to \infty}
  \lambda^n\big(I_{\xi_n}(t)-I_{\xi_n}(T)\big) A_{\xi_n}^{-1}(T)y
\nonumber\\
=\!
  \lambda^0\big(I_0(t)\!-\!I_0(T)\big) A_0^{-1}(T)y\geqref{mmm} \!\!\lambda^0\!\!\int_{T}^{t}\!
    \frac{\partial g}{\partial
    x}(\vartheta,x^0(\vartheta),u^0(\vartheta))\frac{ y\,d\vartheta}{r_0(\vartheta)}
\! \geq \! 0\quad\forall
    t>T.\qquad \label{2023}
\end{eqnarray}
 Substituting $T=0$ and passing to the limit as $t\to\infty,$ we obtain
 the lower estimate from \rref{2455}.

If $\lambda^0>0$, and, in addition, there exists the Lebesgue point
$t^*$ with the required property, then for all points $T\leq
t^*,t>t^*,$
sufficiently close to $t^*$, integration on $[T,t]$
yields "$>$"\ instead of
 "$\geq$"\ in the latter inequality  of \rref{2023}.
Since by \rref{mmm} this integrand is nonnegative for all
$t\in\mm{T}$, the same is true for all $T\leq t^*,t>t^*,$ whence we
obtain $
   \psi^0|_{[0,t^*]}\succ_{\fr{C}} 0_\mm{X}$.

Regarding the latter point, note that if we have $\psi(t)\not\succ_{\fr{C}} 0_\mm{X}$
for some $t\in\mm{T},$ then taking $t^*$ from $\la t,\infty\ra$ yields a contradiction. \bo

\begin{remark}
 \label{mon7}
For $\psi^0(0)\succ_{\fr{C}}0$, it is sufficient to find for each
vector $y_i$ from some basis of~$span\, \fr{C}$  Lebesgue point
$t^*_i$ with the property
 $\frac{\partial g}{\partial x}(t^*_i,x^0(t^*_i),u^0(t^*_i))\, y_i>0$.
\end{remark}

 Let the right-hand side of the dynamics equation and
 the integrand of the objective functional be monotonic in $x$.
  This case frequently arises in economical
 applications, and monotonicity simplifies its analysis.
 It seems that the first to note the peculiarities of this case
  and to investigate it were Aseev, Kryazhimskii, and Taras'ev
  in their paper~\cite{kr_as_t}. These were followed
  by papers~\cite{kr_as2,weber};
the most general cases were considered in~\cite{kr_as,kab}.

Fix the cone $\fr{C}\rav \mm{T}^{\dim \mm{X}}$. In this case,
$\fr{C}^{\bot}=\fr{C}$. Replace $\succcurlyeq_{\mm{T}^{\dim
\mm{X}}},$ $\succ_{\mm{T}^{\dim \mm{X}}}$ with $\succcurlyeq, \succ$.
We obtain
\begin{corollary}
\label{mon}
 Assume conditions $\bf{(u_\sigma),(fg)}, (\tau)$ hold.
  Assume that, for all $x\in\mm{X}$ and for a.a. $t\in\mm{T}$,
   the matrix $\frac{\partial f}{\partial x}(t,x,u^0(t))$ is a matrix with nonnegative off-diagonal
  entries, and $\frac{\partial g}{\partial x}(t,x,u^0(t))$ is a
  nonnegative covector, i.e.,
   there exists a number $d(t,x)\in\mm{R}$  such that the following relation holds:
\begin{equation}
   \label{2433}
 \frac{\partial g}{\partial x}(t,x,u^0(t))\succcurlyeq
0_\mm{X},\quad
  \frac{\partial f}{\partial x}(t,x,u^0(t))\succcurlyeq
  d(t,x)1_\mm{L}.
\end{equation}
Then, there exists a $\tau${-vanishing} multiplier
 $(\lambda^0,\psi^0)\in\Lambda$,
  {and for every such multiplier we have} 
  $\psi^0\succcurlyeq 0_\mm{X},$ and
\begin{equation}
 \label{2455_}
     \lambda^0 \limsup_{t\to\infty,\xi\to 0_\mm{X}} I_\xi(t)
    \succcurlyeq \psi^0(0) \succcurlyeq \lambda^0
   \lim_{t\to\infty} I_0(t) \succcurlyeq 0_\mm{X}
\end{equation}
 holds, and all limits in \rref{2455_} are correctly defined and finite.

If $\lambda^0>0$ (for example, under \rref{1846}) and there exists a Lebesgue point $t^*\in\mm{T}$ for the control~$u^0$ such that
$$  \frac{\partial g}{\partial x}(t^*,x^0(t^*),u^0(t^*))\succ 0_\mm{L},$$
we have
 $
   \psi^0|_{[0,t^*]}\succ 0_\mm{X}$; in particular, $\psi^0(0)\succ 0_\mm{X}.$
\end{corollary}
\begin{remark}
\label{mon12} Assume that under conditions of 
{Corollary}~\ref{mon} we can choose $d(t,x)\equiv 0$, and the
integral
 $$\int_{0}^{t}
    \frac{\partial g}{\partial
    x}(\vartheta,x^0(\vartheta),u^0(\vartheta))\, d\vartheta $$
unboundedly increases as $t\to\infty$; then, all $\tau$-vanishing
solutions are degenerate.
\end{remark}
 Indeed, under $d(t,x)\equiv 0$, we can assume $r_0\equiv 1$; then, in the case $\lambda^0>0$, \rref{2023} would, for $T=0$, imply the boundedness of this integral.

 Note that in \cite[Theorem 1]{kr_as_t}, \cite[Theorem
10.1]{kr_as},\cite[Theorem 5]{kab}
the estimate $\psi\succcurlyeq 0_\mm{X}$  is
made for problems
 \begin{equation}
 \label{eco}
\dot{x}=f(x,u),u\in U,x(0)=x_0,\quad\int_0^\infty e^{-\rho
t}g\big(x(t),u(t)\big)dt\to max.
 \end{equation}
The most general case is examined in \cite[Theorem 5]{kab}; namely, a variant of
 {Corollary} \ref{mon} is stated:
if \rref{2433} is satisfied for all
  $t\in \mm{T},u\in U(t), x\in\mm{X}$ (see~ \cite[(A8)]{kab}), then $\psi\succcurlyeq 0_\mm{X}$, and estimate
  \rref{2455_} holds (see~ \cite[(5.5)]{kab}); the conditions, under which $\psi\succ 0_\mm{X}$ holds in addition to the above, are also specified. The explicit form of estimate \rref{2455_} under the very strong conditions on $f$ and $g$ is
  also specified in \cite[(23)--(26)]{weber}.

Let us also remark that in all papers mentioned, the nondegeneracy of the problem was not assumed (and was not directly reduced to inequality
  \rref{1846}), it had to be proved; for example in \cite[Theorem 5]{kab}, it is demonstrated with the aid of the stationarity condition from additional proposition \cite[(A7)]{kab}:
on any admissible trajectory, there are some $(t,u)$, for which  $f(x(t),u)\succ
  0_\mm{X}$.


\section{Explicit form of $\tau$-vanishing shadow price} \label{sec:6}
\ee

Previously, we examined two transversality conditions \rref{trans} and \rref{partlim}; let us now consider the two conditions
 \begin{subequations}
\begin{eqnarray}
   \label{lim}
       \displaystyle \lim_{t\to\infty}
       ||\psi^0(t)A_0(t)||_{\mm{X}}=0,\\
   \label{partlim_1}
       \liminf_{n\to\infty}
       ||\psi^0(\tau_n)A_0(\tau_n)||_{\mm{X}}=0.
   \end{eqnarray}
          \end{subequations}
\begin{lemma}
\label{lem1}
For each solution $(x^0,u^0,\lambda^0,\psi^0)\in\fr{Y}$,
 the transversality condition~\rref{partlim_1} holds iff $\psi^0(0)$ is a partial limit of the sequence $(\lambda^0 I_0(\tau_n))_{n\in\mm{N}}.$
  \end{lemma}
Passing to the limit in $\lambda^0 I_0(\tau_n)=\lambda^0
(I_0(\tau_n)-I_0(0))=\psi^0(0)A_0(0)-\psi^0(\tau_n)A_0(\tau_n),$ we
obtain what was required; $\lambda^0\neq 0$ by virtue of \rref{dob}.
\begin{lemma}
\label{lem3} If a nontrivial Lagrange multiplier $(1,\psi^0)$
associated with $(x^0,u^0)$  satisfies
 the transversality condition~\rref{partlim_1}, then this multiplier  is
 $\tau$-vanishing.
  \end{lemma}
Indeed, there exists $\tau'\subset\tau$, for which
  $\psi^0(\tau'_n)A_0(\tau'_n)\to 0_\mm{X}.$ Then
  $\psi^0(0)-I_0(\tau'_n)=\psi^0(\tau'_n)A_0(\tau'_n)\to 0_\mm{X},$
  and $I_0(\tau'_n)\to \psi^0(0)$.
Set $\psi_n(t)\rav(I_0(\tau'_n)-I_0(t))A^{-1}_0(t)$.
 Then $\psi_n(\tau'_n)=0_\mm{X}$,
 $\psi^0(0)-\psi_n(0)=\psi^0(0)-I_0(\tau'_n)\ravref{cauchy}\psi^0(\tau'_n)A_0(\tau'_n)\to 0_\mm{X}.$
The proof is completed by virtue of the uniform on each compact
convergence
  $\psi_n\to\psi^0.$

\subsection{Uniformity in initial conditions.} \label{sec:61}

  \begin{theorem}
\label{khlopin_aff_}
 Assume that conditions $\bf{(u_\sigma),(fg)}, (\tau)$ hold.
Let one of the two conditions
\begin{subequations}
   \begin{eqnarray}
   \label{khlopin_omega2}
  either\  \exists I_*&\rav&\lim_{n\to\infty,\xi\to 0_\mm{X}}
I_{\xi}(\tau_n)\in\mm{X};\  \!\!\!\!\!\!\!\!\!\!\!\!\  \\
   \label{khlopin_omega1}
   or\  \exists\iota_*&\rav& \lim_{n\to\infty, \xi\to 0_\mm{X}}
     \frac{I_\xi(\tau_n)}{||I_\xi(\tau_n)||_{\mm{X}}}\in\mm{X},
     \ \!\!\!\!\!\!\!\
     \lim_{n\to\infty, \xi\to 0_\mm{X}}
  ||I_\xi(\tau_n)||_{\mm{X}}=\infty
   \end{eqnarray}
\end{subequations}
hold.

Then, there exists a $\tau${-vanishing} multiplier
 $(\lambda^0,\psi^0)\in\Lambda$.
  Moreover, this multiplier satisfies for all $ T\in\mm{T}$ the corresponding formula of
\begin{subequations}
  \begin{eqnarray}
   \label{klass_}
   \lambda^0= 1,&\ &
 \psi^0(T)\rav
 \Big(I_*-\int_0^T
 \frac{\partial g}{\partial x}(t,x^0(t),u^0(t))\,A_0(t)\,dt\Big)
 A_0^{-1}(T);\\
   \label{khlopin_omega3}
\lambda^0=0,&\ &  \psi^0(T)\rav\iota_* A_0^{-1}(T).
   \end{eqnarray}
\begin{corollary}
\label{aff_}
 Assume conditions $\bf{(u_\sigma),(fg)},(\tau)$ hold.
Let the limit
$$   \lim_{t\to\infty,\xi\to 0_\mm{X}} I_{\xi}(t)=\int_0^\infty
 \frac{\partial g}{\partial x}(t,x^0(t),u^0(t))   A_0(t) dt$$
be well-defined and finite.

 Then, the pair $(x^0,u^0)$ is normal and there exists a unique
 $\tau$-vanishing multiplier $(\lambda^0,\psi^0)\in\Lambda$. Moreover,
  for every solution $(x^0,u^0,\lambda^0,\psi^0)$
  of core relations of the Maximum  Principle \rref{sys_x}--\rref{maxH} and \rref{dobb}, the
   following conditions are equivalent:
 \begin{description}
\item[1)]  its Lagrange multiplier $(\lambda^0,\psi^0)$ is $\tau$-vanishing;
 \item[2)] the transversality condition~\rref{partlim_1} holds;
 \item[3)]  the transversality condition~\rref{lim} holds;
\end{description}
  4)

\vspace{-33.5pt plus 0pt minus 0pt}
  \begin{equation}
  \label{klass} 
   \lambda^0\rav 1,\quad
 \psi^0(T)\rav
 \int_T^\infty
 \frac{\partial g}{\partial x}(t,x^0(t),u^0(t))  A_0(t) \,dt\, {A_0^{-1}(T)}\qquad\forall T\in\mm{T}.
\end{equation}
\end{corollary}
\end{subequations}
  \end{theorem}

Case $(b)$ of {Theorem} \ref{khlopin_aff_} is shown in Corollary
\ref{2353}, case $(a)$ will be proved below together with
{Proposition} \ref{aff}.

In contrast with $(a)$, case~$(b)$ expresses the $\tau$-vanishing Lagrange multiplier
 of a degenerate problem; the author has no knowledge of similar results. Together, these two cases allow to solve problem~\rref{sys}--\rref{opt} through relations of the Maximum Principle regardless of its degeneracy (see, for example, Example \ref{linlin}).


The alternative \rref{khlopin_omega2}$\Rightarrow$\rref{klass_} {\it
vs} \rref{khlopin_omega1}$\Rightarrow$\rref{khlopin_omega3}  is
sufficiently convenient. The need for existence of the limit as
$n\to\infty$ in one of relations
\rref{khlopin_omega2},\rref{khlopin_omega1} can always be satisfied
if we consider a subsequence. However, Example \ref{arnarn} shows
that a unique $\tau$-vanishing multiplier  does not necessarily
satisfy~\rref{partlim_1}, even for {normal} problems.

 Then, the limit in \rref{khlopin_omega2} (or \rref{khlopin_omega1})
should exist not only for $\xi= 0_\mm{X}$, but also as $\xi\to
0_\mm{X}$. In some cases it is provided outright, for example, if the
functions $f$ and $g$ are linear by $x$ (see~
Example~\ref{linlin}), or (see~Example~\ref{avav}) by
the following remark:
 \begin{corollary}
\label{khlopin_aff_5}
Assumptions of Theorem~\ref{khlopin_aff_} hold for a subsequence
  $\tau'\subset\tau$ if one of  the assumptions
  either the functions $f,g$  are linear with respect to $x$,
   \begin{eqnarray*}
 \textrm{or} &\ &   \lim_{n\to\infty, \xi\to 0_\mm{X}}
  \big(I_\xi(\tau_n)-I_0(\tau_n)\big)=0_\mm{X},  
    \\
   \textrm{or} &\ &
   \lim_{n\to\infty, \xi\to 0_\mm{X}}
  \frac{I_\xi(\tau_n)-I_0(\tau_n)}{||I_0(\tau_n)||_\mm{X}}=0_\mm{X},
%
   \end{eqnarray*}
   is satisfied.
  \end{corollary}

\medskip

 Let us finish the proof of Theorem \ref{khlopin_aff_}.
Substituting $T=0$ into \rref{klass_} yields $I_*=\psi^0(0)$; then,
Lemma~\ref{lem1} implies
\begin{lemma}
\label{lem2} A solution $(x^0,\psi^0)$ of
\rref{sys_x}--\rref{sys_psi} given by formula \rref{klass_}
satisfies~\rref{partlim_1} iff~$I_*$ is a partial limit of
the sequence $(\lambda^0 I_0(\tau_n))_{n\in\mm{N}},$
  \end{lemma}
 \begin{proposition}
\label{aff}
 Assume conditions $\bf{(u_\sigma),(fg)},(\tau)$ hold.
Let the map~$I_0$ be bounded and let
 $$\lim_{\xi\to 0_\mm{X}} ||I_{\xi}-I_0||_C=0.$$
 Then, the pair $(x^0,u^0)$ is normal and
 \begin{description}
\item[1)] there exists a $\tau$-vanishing
 multiplier $(1,\psi^0)\in\Lambda$   such that
  transversality condition~\rref{partlim_1} holds;
 \item[2)] a Lagrange multiplier $(\lambda^0,\psi^0)$ associated with $(x^0,u^0)$
 is $\tau$-vanishing iff
the transversality condition~\rref{partlim_1} holds.
\item[3)] a limit point $I_*\in\mm{X}$
  of the sequence   $(I_0(\tau_n))_{n\in\mm{N}}$ corresponds
  to each $\tau$-vanishing multiplier $(\lambda^0,\psi^0)\in\Lambda$, and
  a $\tau$-vanishing multiplier $(\lambda^0,\psi^0)\in\Lambda$ corresponds
  to each limit point $I_*\in\mm{X}$
  of the sequence   $(I_0(\tau_n))_{n\in\mm{N}}$.
  This bijection is given by~\rref{klass_}.
\end{description}
\end{proposition}
\doc
By Theorem \ref{2190}, a $\tau$-vanishing multiplier exists;  by
Remark \ref{2351111}, any $\tau$-vanishing multiplier
$(\lambda^0,\psi^0)$ satisfies $\lambda^0>0$; moreover, by \rref{dobb},
if $(\lambda^0,\psi^0)\in\Lambda$, then $\lambda^0=1$. Now, by
\rref{deimos_1}, we have
$$\psi^0(0)=\lim_{n\to\infty} \lambda^n I_{\xi_n}(\tau'_n)=
  \lambda^0 \lim_{n\to\infty,\xi\to 0_\mm{X}}  I_{\xi}(\tau'_n)=\lambda^0 I_*,$$
and from Lemmas \ref{lem1} and \ref{lem2},  we obtain
\rref{partlim_1} and \rref{klass_}. The inverse is true by virtue of
Lemma \ref{lem3}.
 \bo

\subsection{Uniformity by control}

Formulations of the preceding section can be expressed in another form. By varying, instead of the
initial point~$\xi$, the control $u$ near~$u^0$, we pass from
$x_\xi,A_\xi, I_\xi$ to $x^u,A^u,I^u.$

 Fix pair $(p,\nu)\in\la 0,\infty\ra \times B_{loc}(\mm{T},\mm{R}_{>0})$.
As in Remark~\ref{2351111}, we have
\begin{corollary}
 \label{2350_}
 Assume conditions $\bf{(u),(fg)},(\tau)$.
If for the control $u^0$ and some subsequence $\tau'\subset\tau$ we have
 $$
 \limsup_{n\to\infty,
 \varrho\big(\eta,u^0;\td{\ct{L}}^{p}_{\nu}([0,\tau'_n],\mm{U})\big)\to 0}
   \bigg|\bigg| \int_0^{\tau'_n}
 \frac{\partial g}{\partial
 x}(t,x^\eta(t),u(t))\,A^\eta(t)\,dt\, \bigg|\bigg|_\mm{X}<\infty,$$
then the pair $(x^0,u^0)$ is normal; there exists
  a $\tau$-vanishing multiplier $(1,\psi)\in\Lambda$.
 \end{corollary}
\doc  By Remark~\ref{sas}, there exist a $\tau'$-vanishing multiplier $(\lambda^0,\psi^{0})$
  and sequences $\tau''\subset\tau',(x^n,\eta^n,\lambda^n,\psi^n)_{n\in\mm{N}}$
  such that
 Remark~\ref{1504} and $\varrho\big(\eta^n,u^0;\td{\ct{L}}^{p}_{\nu}(\mm{T},\mm{U})\big)
 \to 0$ hold.
 Then,
 $\varrho\big(\eta^n,u^0;\td{\ct{L}}^{p}_{\nu}([0,\tau''_n],\mm{U})\big)\to
 0$;
 therefore, $(I^{\eta^n}(\tau'_n))_{n\in\mm{N}}$ is bounded by the assumption of the corollary.
 But
 $\lambda^n I^{\eta^n}(\tau''_n)\to\psi^0(0)$, thus $\lambda^0>0.$
 Now $(1,\psi^{0}/\lambda^0)$ is a $\tau$-vanishing multiplier.\bo

\begin{corollary}
\label{s3_d}
 Assume conditions $\bf{(u),(fg)},(\tau)$ hold.
Let~$I_0$ be bounded and let
   $$
   ||I_0- I^{\eta}||_{C([0,\tau_n],\mm{X})}\to 0\quad \mbox{as}\quad n\to\infty, \varrho\big(\eta,u^0;\td{\ct{L}}^{p}_{\nu}([0,\tau_n],\mm{U})\to 0.$$

Then, the pair $(x^0,u^0)$ is normal, and 
 \begin{description}
\item[1)] 
a $\tau$-vanishing  multiplier
 $(1,\psi^0)\in\Lambda$
 corresponds to each partial limit $I_*\in\mm{X}$ of the sequence
 $(I_0(\tau_n))_{n\in\mm{N}}$ by formula \rref{klass_};
\item[2)] all such multipliers satisfy
transversality condition~\rref{partlim_1}.
\end{description}
 \end{corollary}
\doc
Let~$I_*$ be the limit of $(I_0(\tau'_n))_{n\in\mm{N}}$  for certain
  $\tau'\subset\tau.$ Then, by Corollary \ref{2350_},
  there exists a $\tau'$-vanishing multiplier
 $(1,\psi^0)$ such that $\displaystyle\psi^0(0)=\lim_{n\to\infty}I^{\eta^n}(\tau''_n)$
 for some $\tau''\subset\tau'$.
  By the assumption of the corollary this, limit corresponds with
  $I_*$, i.e., $\psi^0(0)=\lambda^0 I_*$.
   But this, by Lemma~\ref{lem1}, is equivalent to \rref{partlim_1}.
   Substituting $\psi^0(0)=\lambda^0I_*$ into \rref{cauchy}, we obtain \rref{klass_}.
    \bo
Repeating the proof of {Corollary} \ref{2350_}, but, this time, using
\rref{1782}, we have
\begin{corollary}
\label{s3_d_d}
 Assume conditions $\bf{(u),(fg)},(\tau)$ hold.
Let for some $\iota_*\in\mm{X}$ there be
   $$%
   \frac{I^{\eta}(\tau_n)}{||I^{\eta}(\tau_n)||_\mm{X}}\to \iota_*,\
      ||I^{\eta}||_{C([0,\tau_n],\mm{X})}\to\infty\quad \mbox{as}\quad n\to\infty, \varrho\big(\eta,u^0;\td{\ct{L}}^{p}_{\nu}([0,\tau_n],\mm{U})\to 0.$$

Then, for the pair $(x^0,u^0)$, there exists a  degenerate
$\tau$-vanishing multiplier
 $(0,\psi^0)$  such that
 condition \rref{partlim_1} and formula \rref{khlopin_omega3}
  hold.
  \end{corollary}

\subsection{Conditions guaranteeing convergence to $I_*$.}
\label{ssect_dom}

Let us consider the conditions on the system that are both sufficiently
 easy to check and sufficient to make use of Corollary~\ref{aff_}.
\begin{proposition}
\label{s3_dd_}
 Assume conditions $\bf{(u_\sigma),(fg)},(\tau)$ hold.
For certain measurable functions
    $F\in \ct{L}^1_{loc}(\mm{T},\mm{L})$, $G\in \ct{L}^1_{loc}(\mm{T},\mm{X}),$ a summable function $\omega\in
   \fr{L}^1(\mm{T},\mm{T})$, let
\begin{subequations}
\begin{equation}
   \label{FG}
   G(t)\succcurlyeq\frac{\partial g}{\partial x}(t,x,u^0(t))\succcurlyeq
   -G(t),\quad
  F(t)\succcurlyeq\frac{\partial f}{\partial
  x}(t,x,u^0(t))\succcurlyeq-F(t),
\end{equation}
\begin{equation}
   \label{Betan}
  ||G(t)B_*(t)||_\mm{X}\leq\omega(t)
\end{equation}
for all
    $(t,x)\in\mm{T}\times\mm{X}$, where $B_*$ is a matrix solution of
\begin{equation}
   \label{Beta}
 \dot{B}_*(t)=F(t)\,B_*(t),\qquad B_*(0)=1_\mm{L}\qquad\forall\,a.a.\, t\in\mm{T}.
\end{equation}
\end{subequations}

Then, the result of Corollary~\ref{aff_} holds.
 \end{proposition}
\doc
For each $B=(b_{ij})_{i,j\in\overline{1,m}}\in\mm{L},$
  $C=(c_{i})_{i\in\overline{1,m}}\in\mm{X}$,
let us introduce $$B^\sharp
\rav(|b_{ij}|)_{i,j\in\overline{1,m}}\in\mm{L},\
  C^\sharp
  \rav(|c_{i}|)_{i\in\overline{1,m}}\in\mm{X}.$$
It is easy to see that $B^\sharp\succcurlyeq 0_\mm{L}$,
$C^\sharp\succcurlyeq 0_\mm{X}$, $B^\sharp\succcurlyeq B \succcurlyeq
-B^\sharp$, $C^\sharp\succcurlyeq C \succcurlyeq
-C^\sharp$.
Moreover,
  $C^\sharp B^\sharp\succcurlyeq CB \succcurlyeq - C^\sharp B^\sharp$
for all $B\in\mm{L},C\in\mm{X}.$

Denote by $F_\xi(t)$ the matrix
   $\frac{\partial f}{\partial x}(t,x_\xi(t),u^0(t))$ for all $t\in\mm{T}$. Now, for all $\xi\in\mm{X}$, we have
\begin{eqnarray*}
  F(t)\succcurlyeq F^\sharp_\xi(t) \succcurlyeq F_\xi(t)\succcurlyeq
  -F^\sharp_\xi(t)\succcurlyeq -F(t)\qquad \forall\ a.a.\ t\in\mm{T};
\end{eqnarray*}
comparing the right-hand sides and the initial
conditions of equations \rref{Beta}, \rref{sys_la}, and equation
$$
 \dot{B}_\xi(t)=F^\sharp_\xi(t)
 B_\xi(t),\qquad B_\xi(0)=1_\mm{L},$$
for its solution $B_\xi$ by the comparison theorem, we obtain
  $$B_*(t)\succcurlyeq B_\xi(t)\succcurlyeq A_\xi(t)\succcurlyeq -B_\xi(t)\succcurlyeq -B_*(t)
  \qquad \forall
  \ a.a.\ t\in\mm{T};$$
 in particular, $B_*(t)\succcurlyeq A^\sharp_\xi(t).$ Now, we have
  $G(t) B_*(t)  \succcurlyeq \big(\frac{\partial g}{\partial
  x}(t,x_\xi(t),u^0(t))\big)^\sharp A^\sharp_\xi(t)
  \succcurlyeq \big(\dot{I}_\xi(t)\big)^\sharp,
  $
whence we obtain $G(t) B_*(t)\succcurlyeq\dot{I}_\xi(t)\succcurlyeq
-G(t) B_*(t)$,
  $||\dot{I}_\xi(t)||_\mm{X}\leq ||G(t) B_*(t)||_\mm{X}\leq \omega(t)$
for all $\xi\in\epsi_0\mm{D}$, for almost all $t\in\mm{T}$. We have
 \begin{eqnarray*}
 ||I_\xi||_C\leq||I_\xi||_{C([0,T],\mm{X})}+\int_{T}^{\infty}\omega(t)\,dt, \\
 ||I_\xi-I_{0}||_C\leq||I_\xi-I_{0}||_{C([0,T],\mm{X})}+2\int_{T}^{\infty}\omega(t)\,dt.
\end{eqnarray*}
For each $\epsi>0$, it is possible to find $T\in\mm{T},$ for which
the second summands do not exceed $\epsi$, and yet
  $I_\xi|_{[0,T]}\to I_{0}|_{[0,T]}$ for $\xi\to 0_\mm{X}.$
  Then all conditions of {Corollary}~\ref{aff_} hold. \bo

\begin{remark}
\label{s3_ddd} The first condition of  \rref{FG} of Proposition
\ref{s3_dd_} could be formally weakened down to
\begin{equation*}
  F(t)+m(t)1_\mm{L}\succcurlyeq\frac{\partial f}{\partial
  x}(t,x,u^0(t))\succcurlyeq-F(t)-m(t)1_\mm{L},
\end{equation*}
for some summable function $m\in
   \fr{L}^1(\mm{T},\mm{T}).$
   \end{remark}
Indeed, consider a number
  $R=e^{\int_0^\infty m(\theta)d\theta}\in\mm{T},$
a summable function
   $\omega_1\rav R\omega$, and a matrix function
   $F_1\rav F +m 1_\mm{L}$. Now,
  $B_{1}(t)\rav e^{\int_{[0,t]}m(\theta)d\theta} B_*(t)$
solves the equation $\dot{B}_1=F_1 {B}_1,{B}_1(0)=1_\mm{L}$
and $$||G(t)B_1(t)||_\mm{X}=e^{\int_{[0,t]}m(\theta)d\theta}||G(t)B_*(t)||_\mm{X}\leq
  e^{\int_{[0,t]}m(\theta)d\theta}\omega(t)\leq R \omega(t)=\omega_1(t).$$
Thus, under conditions of the remark, all propositions of
{Proposition} \ref{s3_dd_} hold for $F_1,\omega_1$
 in the place of $F,\omega.$

Note that conditions of {Proposition} \ref{s3_dd_} (taking into
account {Remark} \ref{s3_ddd}) for a smooth control problem without
phase restrictions are weaker than conditions \cite[(C1)-(C3)]{norv}.
To be more precise, condition \cite[(C1)]{norv} is exactly condition
${\bf{(u)}}$, and
 \cite[(C2)]{norv} is exactly \rref{Betan}. Condition \cite[(C3)]{norv} requires
   $||G(t)B_*(t)B_*^{-1}(\theta)||_\mm{X}\leq\omega(t)$
for all $t\in\mm{T},\theta\in[0,t],$
while condition \rref{FG} requires this only for
  $t\in\mm{T},\theta=0.$
In particular, in \cite[Example 16.1]{kr_as}, conditions of
  \cite[Theorem 12.1]{kr_as} and Proposition~\ref{s3_dd_} hold if $\rho>0,$
   and conditions \cite[(C1)-(C3)]{norv} only hold if
  $\rho>1.$

\begin{corollary}
\label{s3_dd}
 Assume conditions $\bf{(u),(fg)},(\tau)$ hold.
For a summable function $\omega\in \fr{L}^1(\mm{T},\mm{T})$ for all
$u\in\fr{U}$, let
  \begin{equation}
   \label{exp1}
  \Big|\Big|\frac{\partial g}{\partial x}(t,x^u(t),u(t))\, A^u(t)\Big|\Big|_\mm{X}\leq
  \omega(t).
\end{equation}

Then, the pair $(x^0,u^0)$ is normal and Corollary~\ref{aff_} holds
with exception of uniqueness of the $\tau$-vanishing multipliers;
specifically,
 \begin{description}
 \item[1)] exactly one $\tau$-vanishing multiplier
 satisfies \rref{dobb} and \rref{partlim_1};
  \item[2)] exactly one $\tau$-vanishing  multiplier
 satisfies \rref{dobb} and~\rref{lim};
  \item[3)] actually, it is the $\tau$-vanishing multiplier $(1,\psi^0)\in\Lambda$;
   and this multiplier could be obtained by formula
  \rref{klass}.
\end{description}
 \end{corollary}
\doc
Note that \rref{exp1} holds not only for all $u\in\fr{U},$ but also
for all $\eta\in \td{\fr{U}}$; then, for all $T\in\mm{T}$, we have
\begin{eqnarray*}
 ||I^\eta||_C\leqref{exp1}||I^\eta||_{C([0,T],\mm{X})}+\int_{T}^{\infty}\omega(t)\,dt, \\
 ||I^\eta-I_0||_C\leqref{exp1}||I^\eta-I_0||_{C([0,T],\mm{X})}+2\int_{T}^{\infty}\omega(t)\,dt.
\end{eqnarray*}
For each $\epsi>0$ there exists a $T\in\mm{T}$ such that the second
summands do not exceed $\epsi/2.$ Let us construct the
$\tau$-vanishing multiplier $(\lambda^0,\psi^0)\in\Lambda$  by a
limit of sequences from
   Remark~\ref{1504}, but
Proposition~\ref{spectr} implies
  $I^\eta|_{[0,T]}\to I_0|_{[0,T]}$ for $\eta\to\td{u}^0.$ Hence,
   $||I^{\eta^n}-I_0||_C\to 0$ and $I_0$ is bounded.
%
   Since $\psi^n(\tau'_n)=0_\mm{X}$,
    we know that \rref{partlim_1} holds  for $\psi^0$.

From \rref{exp1} for $u=u^0$ we see that for any unboundedly
increasing sequence of times $\upsilon$, the sequence
  $(I_0(\upsilon_n))_{n\in\mm{N}}$ is fundamental and thus it has the limit point~$I_*$.
   Since this is correct for any unboundedly increasing sequence of times, $I_0(t)\to I_*$
   as $t\to\infty$. Lemma \ref{lem1} yields item 2).
Finally, Lemma \ref{lem2} implies \rref{klass}.\bo

The  formula \eqref{klass} was obtained by Kryazhimskii and Aseev
under easily checked assumptions on growth of functions
  $f,g$ and their derivaties (see  stationary case in~\cite[Theorem~12.1]{kr_as},
  \cite[Theorem~4]{kab} and
  non-stationary case in \cite[Theorem~1]{av}).
   This condition
generalizes (see~\cite[Sect.~16]{kr_as}, \cite{av}) a number of
transversality conditions; in particular, it is more general than the
conditions that were obtained for linear systems in~\cite{aucl}.

 From conditions of
\cite[Theorem~2]{kr_as3},\cite[Theorem~12.1]{kr_as}, and
\cite[Theorem~1]{kr_asD} it follows that for some $\alpha,\beta>0$
and for all admissible controls~$u\in\fr{U}$, all trajectories~$x$,
and all fundamental matrices~$A$, the following inequality
holds:
  \begin{equation}
   \label{exp}
  \Big|\Big|\frac{\partial g}{\partial x}(t,x(t),u(t))\Big|\Big|\, ||A(t)||\leq \beta e^{-\alpha
  t}\qquad
   \forall t\in\mm{T}
\end{equation}
  (see, for example,  \cite[(A5)-(A7)]{kr_as}). This is stronger than the conditions of
   Corollary~\ref{s3_dd}.
In  paper \cite{av}, it was actually assumed that \rref{exp} holds
  for $x=x_\xi, A=A_\xi,u=u^0$ if $\xi$ is sufficiently small.
  This is slightly stronger than the stability condition in Corollary~\ref{aff_}.
  However, it is worth noting that \cite[Theorem~1]{av} uses a more general definition of optimality (the locally weakly overtaking optimality).
In addition, condition \rref{exp} can be verified by calculating the characteristic
Lyapunov exponents of the system of the Maximum Principle,
see~\cite[Sect.~12]{kr_as},\cite[Sect.~3]{kr_as3},\cite[Sect.~5]{av}.

Observe that \rref{exp} are characteristic of economic  problems with
exponentially decreasing discount factor; however, one could consider
other non-subexponential discount
  factors
  (see \cite{ivar1,ivar2,weit1,weit2}). Example \ref{avav} exhibits the solution of a problem with such
   discount  factor.

 For economic problems with decreasing discount factor
(specifically, for \rref{eco}) in \cite[Theorem 4]{kab}, sufficiently
broad conditions for applicability of formula \rref{klass} were obtained.
 It turns out that it is sufficient to connect (see~
\cite[(A4)]{kab} and \rref{2742}) the growth of $I^u$ with the
growth of $J^u$. In contrast with the results of
  \cite{av} or {Corollary} \ref{s3_dd}, the finiteness of the optimal result
  on the optimal trajectory is required, and
it is not guaranteed that the $\tau$-vanishing multiplier is unique.
   Let us transfer this result
  of \cite[Theorem 4]{kab} from case \rref{eco} to
  general non-stationary system \rref{sys}--\rref{opt}.

\begin{corollary}
\label{s3_dd_kab}
 Assume conditions $\bf{(u),(fg)},(\tau)$ hold.
Let there exist the finite limit
$\displaystyle\lim_{n\to\infty}J^{u^0}(\tau_n)$. Let a functions
$\omega_0,\omega_\infty\in C(\mm{T},\mm{T})$ satisfy $\omega_0(0)=0$,
   $\omega_\infty(\tau_n)\to 0$
   as $n\to\infty$.
   For all $u\in\fr{U}$ from some $\ct{L}^p_{\nu}(\mm{T},\mm{U})$-neighborhood $O^p_{\nu}$
     of the control $u^0$ for all $k,n\in\mm{N},k<n$, let there be
  \begin{equation}
   \label{2742}
  \Big|\Big|\int_{\tau_k}^{\tau_n}\frac{\partial g}{\partial x}(t,x^u(t),u(t))\, A^u(t)
  dt\Big|\Big|_\mm{X}\leq
  \omega_\infty(\tau_k)+\omega_0\big(|J^u(\tau_n)-J^{u}(\tau_k)|\big).
\end{equation}
Then, the pair $(x^0,u^0)$ is normal, the limit $\displaystyle
 I_*=\lim_{n\to\infty}I^0(\tau_n)\in\mm{X}$ is well-defined, and
  \begin{description}
\item[1)]
 exactly one multiplier
 satisfies \rref{dobb} and \rref{partlim_1};
  \item[2)] actually, it is the $\tau$-vanishing multiplier $(1,\psi^0)\in\Lambda$;
   and this multiplier could be obtained by formula
  \rref{klass}.
\end{description}
 \end{corollary}
 \doc
There exists a
sequence $(s_k)_{k\in\mm{N}}\downarrow 0$ such that for all
$k,n\in\mm{N},k<n$, we have $|J^{u^0}(\tau_n)-J^{u^0}(\tau_k)|<s_k.$
Substituting $u=u^0$ into \rref{2742} yields the existence of the
finite limit $\displaystyle I_*=\lim_{n\to\infty}I^0(\tau_n).$ Now,
as in the proof of {Corollary} \ref{s3_dd}, we show that there exists
the unique solution from $\fr{Y}$ that satisfies
   \rref{partlim_1} and that for it,
  accurately to a positive factor, the formula \rref{klass_} is correct. It only remains
  to prove that  multiplier  defined by
  \rref{klass_} is $\tau$-vanishing.

By Theorem \ref{2190}, for this problem there exists the
$\tau$-vanishing multiplier
  $(\lambda^0,\psi^0)\in {\Lambda}$ that was constructed by
   the uniform limit of sequences $(x^n,\eta^n,\lambda^n,\psi^n)_{n\in\mm{N}}\in \td{\fr{Y}}^\mm{N}$
    from Remark~\ref{1504}. 
      Passing to the subsequence $\tau'\subset\tau$
   if necessary, we may assume $\eta^n\in cl\,O^p_{\nu}$ for all $n\in\mm{N}$,
    then \tref{2742} hold for each $\eta^n$. The function $\omega_0$ can be
   considered monotonic without loss of generality.
   Then, using the triangle
    inequality twice and, by the inequality $\td{J}^{\eta^n}(\tau_n)-J^{u^0}(\tau_n)\geq 0$,
    for all $k,n\in\mm{N},k<n$, we have the following:
\begin{eqnarray*}
\ \ ||I^{\eta^n}(\tau_n)-I^0(\tau_n)||_\mm{X}-
  ||I^{\eta^n}(\tau_k)-I^0(\tau_k)||_\mm{X}\leq  \\
  \ \ ||I^{\eta^n}(\tau_n)-I^{\eta^n}(\tau_k)||_\mm{X}+||I^{0}(\tau_n)-I^{0}(\tau_k)||_\mm{X}
  \leqtref{2742}
  \end{eqnarray*}
  \vspace{-10mm plus 0mm minus 0mm}
\begin{eqnarray*}
 2\omega_\infty(\tau_k)&+&\omega_0(|J^{u^0}(\tau_n)-J^{u^0}(\tau_k)|)
  +\omega_0(|\td{J}^{\eta^n}(\tau_n)-\td{J}^{\eta^n}(\tau_k)|)
  \leq
\\[2ex]
  2\omega_\infty(\tau_k)&+&\omega_0(|J^{u^0}(\tau_n)-J^{u^0}(\tau_k)|)+\\
  &+&\omega_0\big(\td{J}^{\eta^n}(\tau_n)-J^{u^0}(\tau_n)+
  |\td{J}^{\eta^n}(\tau_k)-J^{u^0}(\tau_k)|
  +|J^{u^0}(\tau_n)-J^{u^0}(\tau_k)|\big)\leqtref{555}\\[2ex]
    2\omega_\infty(\tau_k)&+&2\omega_0\big(\gamma_n^2+
    |\td{J}^{\eta^n}(\tau_k)-J^{u^0}(\tau_k)|+s_k\big).
\end{eqnarray*}
Since $I^\eta,\td{J}^\eta$ converges to $I_0,{J}^{u^0}$ uniformly on
any compact and by definitions of $\omega_0,\omega_\infty$ and
$\gamma_n,r_n,s_k$, passing to the limits first, as $n\to\infty$, and
then, as $k\to\infty,$ we see that
  $$\limsup_{n\to\infty}||I^{\eta^n}(\tau_n)-I^0(\tau_n)||_\mm{X}\leq
  2\omega_\infty(\tau_k)+2\omega_0(s_k)$$ and $I^{\eta^n}(\tau_n)-I^0(\tau_n)\to 0_\mm{X}.$
 Now, by Remark~\ref{1504}, we have $\lambda^n I^{0}(\tau_n)\to \psi^0(0)$.
Since $I^0(\tau_n)\to I_*$, we know that $\lambda^0>0$ and $\lambda^0\psi^0(0)=I_*$
hold. By dividing this $(\lambda^0,\psi^0)$ on $\lambda^0$, we obtain
\rref{klass_}.
 \bo


\section{Examples} \label{sec:8}

\begin{example}
The feature of \cite[Ex. 10.2]{norv} lies in the fact that
transversality condition
  \rref{trans} fails to give any information that could help us in determining the unique Lagrange multiplier.
 Let us show that the definition of a $\tau$-vanishing multiplier allows us to do it.

\label{seisei}
$$\dot{x}=ux,\ x(0)=1,\ u\in[1/2,1],\qquad
      J^u(T)\rav\int_{0}^Txe^{-2t}dt\rai\max.$$
\end{example}
Here, $H=u\psi x+e^{-2t}\lambda x$ and
 $\dot{\psi}=-u\psi-e^{-2t}\lambda.$ Then, $A=x,I^u=J^u$; consider
  $F=1,G=e^{-2t},\omega(t)=e^{-t}$. By Proposition \ref{s3_dd_}, there exists the unique $\tau$-vanishing
multiplier.
 Substituting it into $H$, we obtain $H(x^0(t),t,u^0(t),\lambda^0,\psi^0(t))=u^0
  \lambda^0 (J^{u^0}(\infty)- J^{u^0}(t))+e^{-2t}\lambda x^0(t)$;
now, from \rref{maxH}, we have $u^0\equiv 1,$ $J^{u^0}(+\infty)=1$;
then, $\psi^0(0)=\lambda^0=1$, it is a unique $\tau$-vanishing
multiplier. (Of course, in this example,  the control~$u^0$ is easily found in view of
 the monotonicity of $f,g$ and Corollary \ref{mon}).

\medskip

The alternative  \rref{khlopin_omega2}$\Rightarrow$\rref{klass_}
versus \rref{khlopin_omega1}$\Rightarrow$\rref{khlopin_omega3} allows
us to effectively reduce an optimal problem to the boundary problem
of relations of the Maximum Principle.
 The only obstacle is the uniformity of limits in \rref{khlopin_omega2}
  and \rref{khlopin_omega1}. In some cases, the uniformity of these limits
   is trivial, for example, when the functions $f$ and $g$ are linear by $x$. Thus, such problems are easy to solve.
    Let us demonstrate this by the following example:

\begin{example}
\label{sternstern}
 $$\dot{x}=y,\quad\dot{y}=-x+u,\quad x(0)=1,\ y(0)=0,\quad
 u\in [-1,1],\qquad \int_0^T y dt\rai\max$$
\end{example}
Here, for all  $t,T,s\in\mm{T},\xi\in\mm{X}$, we have
\begin{eqnarray*}
  A_\xi(t)=
 \Bigg(
 \begin{array}{cc}
  \cos t & \sin t \\
 -\sin t & \cos t \\
 \end{array}
 \Bigg),
 I_\xi(T)=
 \left(
 \cos T\! -\! 1,
 \sin T
 \right),\\
 I_0(s)A^{-1}_0(T)
 =
  \big(
 \cos  (s\,-\,T)\,-\,\cos T,
 \sin  (s\,-\,T)\,+\,\sin T
 \big).
\end{eqnarray*}

Now,  because $I_\xi$ is $2\pi$-periodic, for any sequence
$(\tau_n)_{n\in\mm{N}}$
 there exists a $\varsigma\in[0,2\pi]$ and subsequence $\tau'\subset\tau$ such that
  $I_\xi(\tau'_n)\to I_0(\varsigma)$, whence, by Theorem~\ref{khlopin_aff_},
  \begin{eqnarray}
  \psi^0(T)&=&(I_0(\varsigma)-I_0(T))A_0^{-1}(T)=
 (\cos(\varsigma-T)-1,\sin(\varsigma-T));\nonumber\\
  u^0(T)&=&arg \max_{u\in[-1,1]}
  (\cos(\varsigma-T)-1,\sin(\varsigma-T))
  \Bigg(
 \begin{array}{c}
  0\\
  u
 \end{array}
 \Bigg)=arg \max_{u\in[0,1]} \sin(\varsigma-T) u,\ \mbox{i.e.}\nonumber\\
  \label{2658}
&\ &\qquad u^0(T)=sgn \sin(\varsigma-T)\quad \forall \mbox{\ a.a.\ } T\in\mm{T}.
 \end{eqnarray}

Observe that the proposed approach finds, first of all, $\tau$-optimal controls. Indeed, let the sequence~$\tau$ be given. Express each~$\tau_n$ in the form $\tau_n=2\pi k_n+\sigma_n$, where $\sigma_n\in[0,2\pi\ra$.
Substituting each limit point $\varsigma$ of the sequence
$(\sigma_n)_{n\in\mm{N}}$ into \rref{klass_} yields all corresponding $\tau$-vanishing multipliers; moreover, formula \rref{2658} yields all prospective $\tau$-optimal controls.

It is easy to check (see \cite{stern}) that any  control of form \rref{2658} is uniformly  weakly
overtaking optimal, thus  each of them is $\tau$-optimal
for its sequence $\tau.$

Also observe that this example specifies why it is impossible
 to replace transversality condition  \rref{partlim_1} in Proposition \ref{aff}
 with the stronger one \rref{lim}.

\begin{example}
Theorem \ref{khlopin_aff_} allows, in some circumstances,
 to find optimal solutions for degenerate problems in the way
 it is done for nondegenerate. Let us show this.
Consider the modification of the well-known Halkin's example \cite{Halkin}
(see also
 \cite[Ex. 5.1]{Pickenhain}, \cite[Ex.~1]{kr_as_t})
 \label{linlin}
 $$\dot{x}=ux,\quad x(0)=1,\quad \int_{0}^T(1-u)x\, dt\rai \max, \quad u\in [\alpha,\beta]\
 (\alpha\leq \beta).$$
\end{example}
 Let there exist a weakly uniformly overtaking optimal control in
this problem, then, for some sequence $\tau$, this control is
$\tau$-optimal.

 Here, $A_{\xi}(T)=x^{0}(T)$ and $I_{\xi}(T)=
J^{u^0}(T).$  Passing, if necessary, from $\tau$ to its subsequence,
we face one of the three cases:

   A) $J^{u^0}(\tau_n)\to +\infty$.  From
   Theorem~\ref{khlopin_aff_}\,(b)
  $\iota_*=1,$ $\lambda=0,$ $H(T)=u^0,$ $u^0\equiv \beta$; if we substitute this into $J^{u^0}(T)$, we will obtain $0\leq\beta<1$.

 B) $J^{u^0}(\tau_n)\to -\infty$; similarly, we have
 $u^0\equiv \alpha>1.$

 C) $J^{u^0}(\tau_n)\to I_*\in\mm{R}$. Here, by Theorem~\ref{khlopin_aff_}\,  (a), from \rref{khlopin_omega2} follows \rref{klass_}.
Consider $R(t)\rav I_*-J^{u^0}(t)-x^0(t)e^{\rho
   t}$. Now we have $H(t)=R(t)u-x$, and $u^0(t)$
   is defined by the sign of $R(t).$ Since $\dot{R}(t)=-x(t)<0$,
    there is at most one switching point.

Note that $u(t)=\gamma$ for all $t>T$, and for some
   $T\in\mm{T},\gamma\in[\alpha,\beta]$. The boundedness of $I_*-J^{u^0}(t)$ provides that either
   $\gamma<0$  or $\gamma=1$.We claim that the sign of $R(t)$ does not change.
   Assume the converse, and let
   there be a switching point $T>0$; then,
 $R(t)<0$ for $t>T$, and $u(t)=\beta=1$, whence $I_*=J^{u^0}(T)$,
  i.e., $x(T)=-R(T)=0$, which is impossible.
Hence, if $R(0)>0$, then $u^0\equiv\beta=\gamma<0;$ else,
$u^0\equiv\alpha=\gamma=1$.

Checking this, we show that
$u^0\equiv\alpha$ for $\alpha\geq 1$ and
 $u^0\equiv\beta$ for $\beta<1$ are indeed $\tau$-optimal (moreover, even uniformly overtaking optimal)
 control in this problem.
  Consequentially, the problem has no $\tau$-optimal
  (and, therefore, no weakly uniformly overtaking optimal) control
   if $\alpha<1\leq\beta$. On the other hand, in case
  $[\alpha,\beta]\rav [0,1]$, the control $u^0\equiv 0$ is
  decision horizon optimal (DH-optimal, see \cite{car1}).
  Therefore, in Theorem \ref{2190}, we could not replace
   the $\tau$-optimality
  (weakly uniformly overtaking optimality, uniformly overtaking optimality)  with  the
  $DH$-optimality (weekly agreeable, agreeable optimality;
  \cite{car1}).

\begin{example}
\label{arnarn} Consider the Arnold's model from \cite{arn}
$$\dot{x}=u,\  x(0)=x_{**},\  u\in [1,2],\ x\in\mm{R},\quad
\frac{\int_{0}^T g(x)\, dt}{T}\rai \max,$$
\end{example}
where profit density, denoted by $g$, is a scalar 1-periodic smooth function
 with a finite number of critical points.

As shown in \cite{arn1,dav}, this problem has  a unique periodic
optimal solution $u^0$, and for certain $g_*\in\mm{T}$, we have
\begin{equation}
\label{arn1}
 (g(x^0(t))<g_*)\Rightarrow (u^0(t)=2)\quad(g(x^0(t))>g_*)\Rightarrow (u^0(t)=1)\quad for
  \, a.a.\, t\in\mm{T}.
\end{equation}
Denote the period of this solution by~$T_0.$

Consider the sequence
$\tau\rav(nT_0)_{n\in\mm{N}}.$ Note that only
the control $u^0$ is $\tau$-optimal  for the problem
\begin{eqnarray}
\dot{x}=u,\  x(0)=x_{**},\  u\in [1,2],\ x\in\mm{R},
\nonumber \\[-1.0ex]
   \label{arn2}
\\
 J(T)={\int_{0}^T g(x(t))\, dt}\rai \max.\nonumber
\end{eqnarray}
Actually, it is possible to prove that this control is at most
 weakly uniformly overtaking and there are no other
weakly uniformly overtaking optimal
 controls in this problem.

Application of Theorem \ref{1} to problem \rref{arn2} yields \rref{arn1}.
Simple reflections on optimality show that
$\min_{x\in[0,1]}g(x)<g^*<\max_{x\in[0,1]}g(x)$, by
 \rref{maxH} we have  $\lambda>0$ for any
 Lagrange multiplier $(\lambda,\psi)$ associated with $(x^0,u^0)$.
  However, no additional conditions on $g_*$ could be obtained
   from the core relations of the Maximum Principle. Let us see if it is possible to do that
   using the approach of this paper.

It is obvious that  $A_\xi\equiv 1_\mm{L}.$ It is also easy to
see that, using the substitution
$\vartheta(t)=x^0(t),t=\vartheta^{-1}(x^0(t))$, we could obtain for
all $T\in\mm{T}$ the following relation:
 \begin{equation*}
\label{sda}
 I_0(T)=\int_{0}^T\frac{dg}{dx}(x^0(t))dt=
      \int_{x^0(0)}^{x^0(T)}\frac{dg(\vartheta)}{d\vartheta}\frac{d\vartheta}{u^0(t)};
       \end{equation*}
now, if $u^0$ is constant on some interval $\la t_2,t_1\ra$, then
 \begin{eqnarray*}
  I_0(t_1)-I_0(t_2)\!=\!\int_{x^0(t_2)}^{x^0(t_1)}\!
  \frac{dg(\vartheta)}{d\vartheta}\!\frac{d\vartheta}{u^0(t)}\!=\!
  \frac{g(x^0(t_1))-g(x^0(t_2))}{u^0(\frac{t_1+t_2}{2})}.
\end{eqnarray*}
 If $t_1,t_2$ are
switching points, then $g(x^0(t_1))=g(x^0(t_2))=g_*,$
$I_0(t_1)=I_0(t_2).$ Since $u^0$ is $T_0$-periodic, this
immediately yields that the functions $x^0,g\circ x^0,I^0$ are also
$T_0$-periodic. 

Observe that the $\tau$-vanishing multiplier $(1,\psi^0)\in{\Lambda}$
exists. Let us show that it does not necessarily satisfy
~\rref{partlim_1} and
  \rref{klass}.
Since $I_0$ is $T_0$-periodic, $I_0(\tau_n)\equiv I_0(0)$, whence
$I_*=0$. If \rref{klass} holds, then, for all $T\in\mm{T}$, we have
$\psi^0(T)\rav -I_0(T)$. Substitution into the Hamiltonian yields
$\ct{H}(t,x^0(t),u,1,\psi^0(t))= -I_0(t)u+g(x^0(t)).$ Now \rref{maxH}
implies that $u^0(t)$ is determined by the sign of $-I_0(t)$, whence
$g_*=g(x^0(0))=g(x_{**})$. But $g_*$  is independent of the choice of
the initial point on the cycle in auxiliary problem~\rref{arn2}.
Therefore, for a.a. $x_{**}$, formula \rref{klass_} is invalid in
this problem. This trivially implies that a $\tau$-vanishing control
does not necessarily satisfy~\rref{partlim_1}, even for normal
problems.

Is it possible to use the formula  \rref{klass_} to find
 $\tau$-vanishing multipliers in this problem? Strange as it sounds,
 yes.

Observe that the notion of $\tau$-vanishing multiplier, as well as the core relations of the  Maximum Principle (see \cite{geom}), is invariant with respect to coordinate transformations. Let us maximize $\overline{J}^u(T)=\ln (1+J^u(T))$ instead of~$J^u(T)$. Consider the problem
\begin{eqnarray}
\dot{x}=u,\  x(0)=x_{**},\ \dot{y}=g(x),\  y(0)=1,\ u\in [1,2],
\nonumber \\[-1.5ex]
   \label{arn3}
\\
\overline{J}(T)=\ln(1+J(T))=\ln y(T) =\int_{0}^T
\frac{g(x(t))}{y(t)}\, dt \rai \max.\nonumber
\end{eqnarray}
Take an arbitrary control~$u^0$ of form~\rref{arn1}, and let its
period be some~$T_0$. It is easily seen that $\displaystyle
\overline{A}_\xi\equiv \Bigg(
 \begin{array}{cc}
 1 & 0 \\
 I_\xi & 1 \\
 \end{array}
 \Bigg),\\$
\begin{eqnarray*}
  \overline{I}_\xi(nT_0)=\int_{0}^{nT_0}
  \frac{1}{y_\xi(t)}\bigg(\frac{dg(x_\xi(t))}{dx},- \frac{g(x_\xi(t))}{y_\xi(t)}\bigg)
  \Bigg(
 \begin{array}{cc}
 1 & 0 \\
 I_\xi & 1 \\
 \end{array}
 \Bigg) \,
 dt=\\
\int_{0}^{nT_0} \frac{1}{y^2_\xi(t)}
  \big({\dot{I}_\xi(t)y_\xi(t)-{I}_\xi(t)\dot{y}_\xi(t)
  },- {\dot{y}_\xi(t)}{}\big)\,dt=\\
  \frac{\big({I}_\xi(t),1\big)}{y_\xi(t)}\bigg|_{t=0}^{t=nT_0}=
    \frac{\big(n{I}_\xi(T_0),1\big)}{n(y_\xi(T_0)-y_\xi(0))+y_\xi(0)}-(0,1)
    \to
    \Big(\frac{{I}_\xi(T_0)}{y_\xi(T_0)-y_\xi(0)},-1\Big).
\end{eqnarray*}
Now, the theorem of continuous dependence on initial conditions implies \rref{khlopin_omega2}. Thus, {Theorem}~\ref{khlopin_aff_} also holds for problem \rref{arn3} for each control~$u^0$ of form \rref{arn1}, and its proper  $\tau$-vanishing Lagrange multiplier is given by formula \rref{klass_}. Thus, formula \rref{klass_}, under proper coordinate transformation, can be used to solve problem~\rref{arn2},\rref{arn3}, although this yields no additional conditions in comparison with the core relations of the Maximum Principle.

Actually, this is rather reasonable since a control of form \rref{arn1} is
  weakly uniformly overtaking optimal for the objective functional
$$\overline{\overline{J}}(T)\rav\ln(1+\ln(1+J(T)))=\int_{0}^T
\frac{g(x)}{y(1+\ln y)}\, dt \rai \max.$$
Therefore, in this problem, it has a $\tau$-vanishing multiplier; since the definition of      $\tau$-vanishing multiplier is invariant, each control of form~\rref{arn1} has such a multiplier in problems \rref{arn2} and  \rref{arn3} too even if the corresponding controls are not        weakly uniformly overtaking optimal in these problems.

\medskip
Let us show the example of reducing an infinite horizon optimal control problem to the boundary problem.
\begin{example}
\label{avav}
In \cite{av}, the following stylized
 microeconomic problem was considered:
 \begin{eqnarray*}
 \dot{x}(t) &=& - \nu x(t) + u(t),\ \  x(0) = K_0,\ u\geq 0;\\
J^u(T) &=& \int_0^T  e^{-dt}\Big[e^{pt} (x(t))^\sigma -
  \frac{b}{2} (u(t))^2\Big]  dt \rai max.
  \end{eqnarray*}
Here,
 $u(t)$ is the investment,
 $\nu\geq 0$ is the depreciation rate, $K_0>0$ is the given initial capital stock,
 $e^{-dt}$ is
the  discount factor ($d\geq 0$), $e^{pt}\geq 0$ is the (exogenous)
factor of technological advancement ($p\geq 0$),
  $bu^2(t)$ ($b > 0$) is the cost of investment $u(t)$, and $\sigma\in\la 0,1]$
defines the production function.
Under the assumption $d+\nu>
 \frac{p}{2-\sigma}$, it is shown that
 there are no optimal solutions for $p>d+\nu$, and, for
   $p<d+\nu$, each locally weakly overtaking control induces
   a solution of the boundary problem (see \cite{av}).
\end{example}

Consider the following objective functional:
\begin{eqnarray*}
J^u(T) = \int_0^T  g(t) (x(t))^\sigma - h(t)
  \frac{b}{2} (u(t))^2  dt \rai max.
\end{eqnarray*}
Here,  $h(t)$ is the discount factor, $g(t)$ is the product of the
discount factor and the factor of technological advancement.

Suppose that there exists a weakly overtaking optimal control
$u^0$. Then, for some sequence $\tau\uparrow \infty$, this solution is $\tau$-optimal. Hence, there exists a $\tau$-vanishing multiplier
$(\lambda^0,\psi^0)\in\Lambda$.

Now, for all $\xi\in\mm{X}$, we have $A_\xi=e^{-\nu t},$
  $$I_\xi(T)=\int_{0}^T g(t) \sigma x_\xi^{\sigma-1}(t)
   e^{-\nu t}\, dt=\sigma\int_{0}^T g(t)e^{-\nu t}   x_\xi^{\sigma-1}(t)
   \, dt.$$

 Note that $x_\xi(t)-x^0(t)=\xi e^{-\nu t}$; now
 we have
\begin{eqnarray*}
   I_\xi(T)-I_0(T)=
\sigma   \int_{0}^T  g(t) e^{-\sigma \nu t}
   \Big[(x^0(t)e^{\nu t}+\xi)^{\sigma- 1}-
   (x^0(t)e^{\nu t})^{\sigma- 1}\Big]\, dt.
\end{eqnarray*}
It is easy see that
$\big|(r+\xi)^{\sigma-1}\!-\!r^{\sigma- 1}\big|\leq (2^{2\!-\!\sigma}-2)|\xi| r^{\sigma-2}\leq
(2^{2-\sigma}\!-\!2) K_0|\xi| r^{\sigma-1}$ if $2|\xi|<K_0\leq r$.
Since the function $x^0(t)e^{\nu t}$ is monotonically increasing, we obtain
\begin{eqnarray*}
|I_\xi(T)\!-\!I_0(T)|\!\leq\! \bigg|\int_{0}^T\!\!  g(t) e^{-\nu t} (x^0)^{\sigma\!-\! 1} dt\bigg|
(2^{2\!-\!\sigma}\!\!-\!2) K_0|\xi|=|I_0(t)|(2^{2\!-\!\sigma}\!\!-\!2)K_0|\xi|
\end{eqnarray*}
 for all $T\in\mm{T},2|\xi|< K_0.$
Now, by Corollary~\ref{khlopin_aff_5},  considering the subsequence if necessary, we have the conclusion of
Theorem
   \ref{khlopin_aff_}.

   We claim that $(I_0(\tau_n))_{n\in\mm{N}}$ is bounded.
   Assume the converse; then, considering the subsequence if necessary, we come to \rref{khlopin_omega1} and \rref{khlopin_omega3}, whence
   $\displaystyle\lim_{\xi\to 0,\ n\to\infty} I_\xi(\tau_n)=\pm\infty,$ now
   $\iota^*=\pm1,\lambda^0=0$ and by \rref{maxH} we have
   $$ u^0(t)= \arg \max_{u\in\mm{R}_{\geq 0}}e^{\nu t}
     I_0(t)(u- \nu x)=\arg \max_{u\in\mm{R}_{\geq 0}}
     I_0(t)u=\pm\infty,
  $$
which is impossible. This contradiction proves the boundedness of
sequence  $(I_0(\tau_n))_{n\in\mm{N}}$.

 Now there exists a finite limit $I_*$  of
$(I_0(\tau'_n))_{n\in\mm{N}}$ for some $\tau'\subset\tau.$
By  Theorem~\ref{khlopin_aff_}, we have~\rref{klass_},
   $\lambda^0=1,\psi^0(T)=(I_*-I_0(t))e^{\nu t},$
\begin{eqnarray*}
   u(t)= \arg \max_{u\in\mm{R}_{\geq 0}}e^{\nu t}
     (I_*-I_0(t))(-\nu x+u)+ g(t) (x^0(t))^\sigma -
  h(t)\frac{b}{2}u^2 =\\
  \arg  \max_{u\in\mm{R}_{\geq 0}}e^{\nu t}(I_*-I_0(t))u-h(t)
  \frac{b}{2}u^2 =\frac{e^{\nu t}}{b h(t)}(I_*-I_0(t))\quad  for\, a.a.\, t\in\mm{T}.
  \end{eqnarray*}
Consider $I(t)\rav I^*-I_0(t)$; differentiating $I(t)$ with respect
to $t$,
 we finally close \rref{sys_x}--\rref{sys_psi} into the boundary problem
\begin{subequations}
\begin{eqnarray}
    \label{3000}
   \dot{x}^0&=&- \nu x^0+\frac{e^{\nu t}}{b h(t)} I,\qquad
   x^0(0)=K_0,\\
    \label{3001}
        \dot{I}&=&-
   {\sigma  g(t)}{e^{-\nu t}} (x^0)^{\sigma- 1}, \\
    \label{3002}
   I(\tau'_n)&\to& 0\mbox{\ as\ }n\to\infty.
\end{eqnarray}
Each $\tau'$-optimal control generates the unique solution of this
problem. For $\sigma=1$ if such solution exists then there exists a
finite limit
 $\displaystyle \lim_{n\to\infty}\int_{0}^{\tau''_n}e^{-\nu t} g(t)dt$ for some $\tau''\subset\tau'.$

Note that to construct this boundary problem we have to know the subsequence $\tau'\subset\tau$. In terms of the initial sequence~$\tau$, it is only possible to claim that, for a solution $(x^0,I)$ of problem \rref{3000}--\rref{3001},
  $0_\mm{X}$ is the partial limit of the sequence $(I(\tau_n))_{n\in\mm{N}}$. If for some functions $g,h$ for some sequence $\tau\uparrow\infty$ there are multiple $\tau$-solutions, then each of them has its own~$I$ and subsequence~$\tau'$. Also note that if we do not know the sequence~$\tau$, then
  to find a
  weakly uniformly overtaking optimal control, we have to solve problem \rref{3000}--\rref{3001} for the boundary condition
   $$\liminf_{t\to\infty}I(t)\leq 0\leq \limsup_{t\to\infty}I(t).$$

Now suppose that $g(t)\geq 0,h(t)>0$ for a.a. $t\in\mm{T}$. Then,
 there exists the common limit~$I_*$ of all sequences $(I_0(\tau_n))_{n\in\mm{N}}$, and for each
 weakly overtaking optimal control~$u^0$  there exists the unique solution of problem
\rref{3000}--\rref{3001} for the boundary condition
\begin{equation}
   \label{3003}
   I(t)\to +0\mbox{ as }t\to\infty.
\end{equation}
 \end{subequations}

It is possible to find the explicit solution of boundary problem
\rref{3000},\rref{3001},\rref{3003} in some specific cases. For
example, let the discount factor equal $\frac{1}{(1+t)^{4/3}}$,
let the factor of technological advancement be equal to 1. For
$$g(t)=h(t)=\frac{1}{(1+t)^{4/3}},\nu=0,\sigma=1/2,b=\frac{3}{8},
K_0=1$$ we have
   $$x^0(t)=(1+t)^{4/3},\ u^0(t)=\frac{4}{3}(1+t)^{1/3},\
   I(t)=\frac{1}{2(1+t)},\ J^{u^0}(t)=(1+t)^{2/3}.$$

     The discount factor $g(t)=\frac{1}{(1+t)^{4/3}}$ here is not
arbitrary, its power
    $\alpha=3,96/2,94\approx 4/3$ was determined by means of statistic analysis in \cite{weit1}.
A thorough discussion of various discount functions and their properties could be found, for example, in
\cite{ivar1,ivar2,weit2}.
   These papers do not generally
  assume the   discount function to be dominated by a decreasing exponential
   function and do not assume its monotonicity.

\section{Appendix}
\label{sec:7} \ee
{\bf The proof of Proposition~\ref{spectr}.}
For the sake of brevity, let us denote $\td\Pi\rav\prod_{n\in\mm{N}}
\td{\fr{U}}_n$, and let us equip it with Tikhonov topology.
 Let  $\td{\Delta}:\td{\fr{U}}\to\td\Pi$  be given by
 $\td{\Delta}(\eta)\rav\big(\td{\pi}_n(\eta)\big)_{n\in\mm{N}}$
for all $\eta\in\td{\fr{U}}.$ It is a homeomorphism by continuity of
the maps~$\td{\pi}_{n}$ and $\td{\pi}_n\circ\td{\Delta}^{-1}$.

Let $n,k\in\mm{N}, (n>k).$ Then, the space $\td{\fr{U}}_n$ is
included in~$\td{\fr{U}}_k$ by the mapping
  $\td{\pi}^{n}_{k}(\eta)\rav\eta|_{[0,k]}$ for all $\eta\in \td{\fr{U}}_n.$
 By $\td{\pi}^{n}_{k}\circ\td{\pi}^{k}_{i}=\td{\pi}^{n}_{i}$ for all
  $n,k,i\in\mm{N}, (n>k>i),$ we have the projective sequence of the topological spaces
  $
  \{
  \td{\fr{U}}_n,\td{\pi}^{n}_{k}\};$
and we can define the inverse limit \cite[III.1.5]{phillvv},
  \cite[2.5.1]{en}. In our notation, we can write it in the form
    $\varprojlim
    \{\td{\fr{U}}_n,\td{\pi}^{n}_{k}\} \rav\td\Delta(\td{\fr{U}})\subset \td\Pi.$
As shown above,~$\td\Delta$ is a homeomorphism; hence, $\td{\fr{U}}$
is homeomorphous to $\td\Delta(\td{\fr{U}})$. Now, by Kurosh
Theorem~\cite[III.1.13]{phillvv}, the inverse limit
$\td\Delta(\td{\fr{U}})$ of compacts $\td{\fr{U}}_n$ is compact, and
 $\td{\fr{U}}$ is a compact too. Similarly, from \cite[4.2.5]{en} and \cite[IV.3.11]{va} it
follows that $\td{\fr{U}}$ is also metrizable.

Repeating the reasonings without $\td{\ }$ or referring
to~\cite[3.4.11]{en} and~\cite[2.5.6]{en} yields
$\displaystyle{\fr{U}}\cong\varprojlim
\{
  {\fr{U}}_n,\pi^{n}_{k}\}\rav \Delta({\fr{U}})\subset \Pi.$

For each $n\in\mm{N}$, let the mapping ${e}_n:\fr{U}_n\to
\td{\fr{U}}_n$ be given by
    ${e}_n(u)(t)\rav(\td\delta\circ u)(t)=\td{\delta}_{u(t)}$ for all
    $t\in[0,n],u\in\fr{U}_n.$ Since for all $n,k\in\mm{N},
  n>k$ it holds that
  ${e}_k \circ {\pi}^n_k = {e}_n,$
we have the projective system $\{{e}_n,{\pi}^n_k\}$. Passing to the
inverse limit, we obtain the mapping
$e_\Delta:\Delta(\fr{U})\to\td\Delta(\td{\fr{U}})$; from
$e_n\circ\pi_n=\td\pi_n\circ\td\delta$ we have
    $e_\Delta\circ\Delta=\td\Delta\circ\td{\delta}$, and from
  $\td{\fr{U}}_n=cl\, e_n(\fr{U}_n)$ (\cite{va}) we have
  $\td\Delta(\td{\fr{U}})=cl\, e_\Delta\big(\Delta(\fr{U})\big)=cl\, (\td\Delta\circ\td{\delta})(\fr{U});$
now, by continuity of~$\td\Delta^{-1}$, we obtain $\td{\fr{U}}= cl\,
\td{\delta}(\fr{U}).$

The mapping $\td{\fr{A}}[\eta]$ is continuous by virtue of, for
example, \cite[Theorem~3.5.6]{tovst1}; the set
$\td{\fr{A}}[\eta](\td{\fr{U}})$ is compact as a continuous image of
a compact. In what follows, is sufficient to use
    $\td{\fr{U}}= cl\, \td{\delta}(\fr{U}).$

    Replacing~$a$ and the compact~$\Xi$
    with the mapping $(f,g)$ and the compact $\{(x_{**},0)\}$,
    we obtain the continuous dependence on~$\eta$
    for the maps $\ph^\eta,\td{J}^\eta$.
\bo
\medskip
\begin{subequations}
{\bf The proof of Proposition~\ref{dop}.} For all $n\in\mm{N}$, let
us consider set
      $$\bar{G}_n\rav\Big\{ (t,y(t))\,\Big|\,\forall y\in\td{\fr{A}}[\td{u}^0],t\in
      {[0,n] } \Big\};$$
      by the theorem of
continuous dependence of solutions on initial conditions
       this set is compact as a continuous image of
a compact~$\Xi.$
        Therefore, on this set, the function $a(t,y,u^0(t))$
       is Lipshitz continuous with respect to~$y$
       for the certain Lipshitz constant
       $L_n\rav L^a_{\bar{G}_n}\in \ct{L}^1_{loc}(\mm{T},\mm{T})$.
       For all $t\in[0,n]$, define
       $M_n(t)\rav\int_0^t L_n(\tau)d\tau$.
       Note that this function is absolutely continuous and monotonically nondecreasing.

Fix $n\in\mm{N}$; for all $t\in [n-1,n\ra, u\in \mm{U}$, let us
consider the number
 $$    R(t,u)\rav \sup_{ y\in \bar{G}_n}
 \big|\big|a(t,y,u)-a(t,y,u^0(t)\big)\big|\big|_E.
     $$
Here, the norm  is continuous with respect to~$y$ and~$u$,
 and~$y$ assumes values from the compact set; now, for every $u\in \mm{U}$
 by~\cite[Theorem 3.7]{select} the supremum reaches the maximum
for the certain function $y_{max}[u]\in
\ct{L}^1([n,n-1\ra,\bar{G}_n)$. Hence, $R(t,u)$ is measurable with
respect to~$t$ for each $u\in\mm{U}$.

Fix a $t\in [n-1,n\ra$; for each sufficiently small neighborhood
$\Upsilon\subset \mm{U}$, by continuity of $a(t,\cdot,\cdot)$ on
compact $\bar{G}_n\times cl\, \Upsilon$, there exists a function
$\omega^t\in C(\mm{T},\mm{T})$ such that $\omega^t(0)=0$ and
\begin{equation}
  \label{1111}
   \Big|\,\big|\big|a(t,y,u')-a(t,y,u^0(t)\big)\big|\big|-\big|\big|a(t,y,u'')-a\big(t,y,u^0(t)\big)\big|\big|\,\Big|<
   \omega^t\big({
   {||u'\!-\!u''||}}\big)
\end{equation}
holds for every $y\in\bar{G}_n, u',u''\in \Upsilon$. Without loss of
generality,
 assume $R(t,u')\leq R(t,u'')$. 
Now, by definition,
  $R(t,u')\geq
  \big|\big|a(t,y,u')-a\big(t,y,u^0(t)\big)\big|\big|$,
   and, substituting $y\rav y_{max}(u'')(t)$ into~\rref{1111}, we obtain
  $0\leq R(t,u'')-R(t,u')\leq
  \omega^t(||u'-u''||);$
  i.e.,~$R$ is continuous   with respect to the variable~$u$ on
each sufficiently small neighborhood $\Upsilon\subset \mm{U}$;
therefore on~$\mm{U}$ too. Thus, the function $R:
[n-1,n\ra\times\mm{U}\to\mm{T}$ is a Carath\'{e}odory function.

Let us note that by considering all $n\in\mm{N}$, we define the
Carath\'{e}odory function~$R$ on the whole $\mm{T}\times\mm{U}$.
Moreover, by construction, $R(t,u^0(t))\equiv 0$. Hence, it is
correct to define $w^0\in (Null)(u^0)$ by the rule
\begin{equation}
\label{nado}
  w^0(t,u)\rav||u-u^0(t)||+e^{M_n(t)}R(t,u) \qquad
  \forall n\in\mm{N},t\in[n-1,n\ra,u\in\mm{U}.
\end{equation}

Consider arbitrary $n\in\mm{N}$, $\vartheta\in [0,n]$, and
$(\vartheta,y^*_1),(\vartheta,y^*_2)\in
        \bar{G}_n$. 
There exist  solutions $y_1,y_2\in \td{\fr{A}}[\td{u}^0]$ of
equation~\rref{1667}, for the initial conditions
$y_i(\vartheta)=y^*_i$. 
Let us introduce functions
        $$r(t)\rav y_1(t)-y_2(t),\quad
                 W_+(t)\rav e^{M_n(t)}||r(t)||_E\qquad \forall t\in [0,n].$$
By Lipshitz continuity of the right-hand side of~\rref{1667} we
obtain $ ||\dot{r}(t)||_E\geq -L_n(t)
       ||r(t)||_E,$ and
       $$\frac{dW^2_+(t)}{dt}=2L_n(t)W^2_+(t)
       +2e^{2M_n(t)}r(t)\dot{r}(t)
       \geq2L_n(t)W^2_+(t)
       -
       2L_n(t)W^2_+(t)
       =0.
       $$
       Thus, the function~$W_+$ is nondecreasing,
       and finally for all $(\vartheta,y^*_1),(\vartheta,y^*_2)\in
        \bar{G}_n$ we have
\begin{equation}
  \label{1037_}
     ||\beg(\vartheta,y^*_1)-\beg(\vartheta,y^*_2)||_E={W_+(0)}\leq {W_+(\vartheta)}=
      e^{M_n(\vartheta)}||y^*_1-y^*_2||_E.
\end{equation}

 Let us now consider
 $\eta\in\td{\fr{U}},$ $y\in\td{\fr{A}}[\eta],$
         $T\in\mm{T}$ such that
          $\varkappa(\vartheta,y(\vartheta))\in\Xi$ for all $\vartheta\in[0,T]$.
         Fix arbitrary $n\in\mm{N}$ and $\vartheta_1,\vartheta_2\in[0,T]\cap[n-1,n\ra,$
        $\vartheta_1<\vartheta_2$
 there exists the solution $y^0\in\td{\fr{A}}[\td{u}^0]$
        such that
        $y^0(\vartheta_1)=y(\vartheta_1)$. By construction of $\bar{G}_n$, we have
        $(t,y(t)),(t,y^0(t))\in \bar{G}_n$ for all $t\in [\vartheta_1,\vartheta_2].$
         Let us also define
        $$r\rav y^0(t)-y(t),\quad W_-(t)\rav e^{-M_n(t)}||r(t)||_E\qquad \forall t\in [\vartheta_1,\vartheta_2].$$
        Then $W_-(\vartheta_1)=0.$
Now,
 \begin{eqnarray*}
  \frac{dW_-^2(t)}{dt}=
       2e^{-2M_n(t)}r(t)\dot{r}(t)-2L_n(t)W_-^2(t)
       =\\
       2e^{-2M_n(t)}{r}(t)\big(\dot{y}^0(t)\!-\!a(t,{y}(t),u^0(t))+
       a(t,{y}(t),u^0(t))\!-\!\dot{y}(t)\big)-2L_n(t)W_-^2(t)
       \leq \\
          2e^{-2M_n(t)}||r(t)||_E
          \int_{U(t)}R(t,u)
          d\eta(t)+2L_n(t)W_-^2(t)-2L_n(t)W_-^2(t)\leq \\
   2e^{-M_n(t)}W_-(t)
          \int_{U(t)}R(t,u)
          d\eta(t)\leq 2e^{-2M_n(t)}W_-(t)\frac{d
          \fr{L}_{w^0}[\eta](t)}{dt}.\ \ \ \
\end{eqnarray*}
     Since function~$W_-$ is nonnegative,
     for a.~a. $t\in \{t\in[\vartheta_1,\vartheta_2]\,|\,W_-(t)\neq 0\}$ we obtain
 \begin{equation*}
     \frac{dW_-(t)}{dt}\leq
     e^{-2M_n(t)}\frac{d \fr{L}_{w^0}[\eta](t)}{dt}\leq e^{-2M_n(\vartheta_1)}\frac{d
     \fr{L}_{w^0}[\eta](t)}{dt}.
\end{equation*}
This inequality is trivial for
     $t\in[\vartheta_1,\vartheta_2], t<\sup \{t\in[\vartheta_1,\vartheta_2]\,|\,W_-(t)=0\};$
whence, integrating  inequality in $t\in [\vartheta_1,\vartheta_2]$,
we obtain
$$W_-(\vartheta_2)=W_-(\vartheta_2)-W_-(\vartheta_1)  \leq
e^{-2M_n(\vartheta_1)}
          \big(\fr{L}_{w^0}[\eta](\vartheta_2)-\fr{L}_{w^0}[\eta](\vartheta_1) \big). $$
But
$\beg(\vartheta_2,y^0(\vartheta_2))=y^0(0)=\beg(\vartheta_1,y^0(\vartheta_1))=\beg(\vartheta_1,y(\vartheta_1)),$
hence, we have
 \begin{eqnarray}
     \!||\beg(\vartheta_2,y(\vartheta_2))-\beg(\vartheta_1,y(\vartheta_1))||_E=
     ||\beg(\vartheta_2,y^0(\vartheta_2))-\beg(\vartheta_2,y(\vartheta_2))||_E
     \leqref{1037_}\nonumber\\
e^{M_n(\vartheta_2)}||y^0(\vartheta_2)-y(\vartheta_2)||_E=
     e^{2M_n(\vartheta_2)}W_-(\vartheta_2)\ \   \leq\ \ \nonumber    \\
e^{2M_n(\vartheta_2)-2M_n(\vartheta_1)}
          \big(\fr{L}_{w^0}[\eta](\vartheta_2)-\fr{L}_{w^0}[\eta](\vartheta_1)
          \big).\ \ \ \
  \label{1649}
\end{eqnarray}

Fix arbitrary $t\in [0,T]$. For each $\epsi>0$ we can split interval
$[0,t\ra$ into the intervals of the form $[\vartheta',\vartheta''\ra$
such that
 $M_n(\vartheta'')-M_n(\vartheta')<\epsi$
 and
 $[\vartheta',\vartheta''\ra\subset[n-1,n\ra$ for the certain $n\in\mm{N}.$
But,~\rref{1649} holds for every interval, i.e.,
$$||\beg(\vartheta'',y(\vartheta''))-\beg(\vartheta',y(\vartheta'))||_E\leq
    e^{2\epsi}
          \big(\fr{L}_{w^0}[\eta](\vartheta'')-\fr{L}_{w^0}[\eta](\vartheta') \big).$$
Summing for all intervals, by $\beg(0,y(0))=y(0)$ and by the triangle
inequality, we obtain
    $||\beg(t,y(t))-y(0)||_E\leq
    e^{2\epsi}
          \fr{L}_{w^0}[\eta](t)$
          for every $t\in [0,T]$. Arbitrariness of $\epsi>0$
           completes the proof of the Proposition \ref{dop}.
 \bo
\end{subequations}
\medskip
\begin{subequations}
{\bf The proof of {Proposition} \ref{1014}.}
 For every $n\in\mm{N}$, let us consider
the problem
$$   J^\eta(\tau_n)-\gamma_n
\fr{L}_w[\eta](\tau_n)=\int_0^{\tau_n}\int_{U(t)}
  g(t,x^\eta(t),u)d\eta(t)dt-\gamma_n
  \fr{L}_w[\eta](\tau_n)
  \to\max.$$
Here, the functional is bounded from above by the number
$J^{u^0}(\tau_n)+\gamma_n^2$, therefore, it has the supremum. Every
summand continuously depends on~$\eta$, which covers the
compact~$\td{\fr{U}}$; therefore, there is an optimal solution for
this problem in~$\td{\fr{U}}$; let us denote one of them by
$\eta^n$, and its trajectory by $x^n$.

  For every $\gamma\in\mm{T}$ let    the function
  $\ct{H}_{\gamma}:\mm{X}\times\mm{T}\times\mm{U}\times\mm{T}\times\mm{X}\to\mm{R}$
 be given by
 $$\displaystyle \ct{H}_{\gamma}(x,t,u,\lambda,\psi)\rav\ct{H}(x,t,u,\lambda,\psi)-
  \gamma w(t,u).$$
Then, by the  Maximum
Principle~\cite[Theorem~5.2.1]{clarke}, there exists
$(\lambda^n,\psi^n)\in\la0,1]\times C([0,n],\mm{X})$ such that
relation~\rref{dob} and
 the transversality condition  $\psi^n(\tau_n)=0$
 hold,
 and
 \begin{eqnarray}
   \label{sys_max_}
   \sup_{p\in
        U(t)}\!\ct{H}_{\gamma_n}\big(x^n(t),t,p,\lambda^n,\psi^n(t)\big)\!&=&\!
           \int_{U(t)}\!\ct{H}_{\gamma_n}\big(x^n(t),t,u,\lambda^n,\psi^n(t)\big)
   d\eta^n(t),\\
\dot{\psi}^n(t)&=& -
\int_{U(t)}\frac{\partial\ct{H}_{\gamma_n}}{\partial
 x}\big(x^n(t),t,u,\lambda^n,\psi^n(t)\big)
   d\eta^n(t) \nonumber
 \end{eqnarray}
    also hold for a.a. $t\in[0,\tau_n]$.

Let us extend the  $(x^n,\eta^n,\lambda^n,\psi^n)$ to
$[\tau_n,\infty\ra$ by the generalized control
$\td{u}^0|_{[\tau_n,\infty\ra}$. 
Let us denote by~$\fr{Z}^n$ the set of $(x,u,\lambda,\psi)$ that
satisfy relations
  \rref{dob},
 $\td{\mbox{\rref{sys_x}}}$--$\td{\mbox{\rref{sys_psi}}}$
a.~e. on~$\mm{T}$, satisfy relation~\rref{sys_max_} a.~e.  on
$[0,\tau_n\ra$, and possess the property
$\td{u}^0|_{[\tau_n,\infty\ra}=
   \eta^n|_{[\tau_n,\infty\ra}$. Now we have $(x^n,\eta^n,\lambda^n,\psi^n)\in\fr{Z}^n$
for every $n\in\mm{N}$.

Let us note that all~$\fr{Z}^n$ are closed and, since these sets are
contained
 in the compact~$\td{\fr{Y}}$, these sets are also compact. Hence, the sequence
  $(x^n,\eta^n,\lambda^n,\psi^n)_{n\in\mm{N}}$ has the limit point
  $(x^{00},\eta^0,\lambda^0,\psi^0)\in\td{\fr{Y}}$.
Considering, if need be, the subsequence,  we can assume that it is
the limit of the sequence itself.

For all $t,\gamma,\lambda\in\mm{T},(x,\psi)\in\mm{X}\times\mm{X}$,
denote by $\ct{P}_{\gamma,\lambda}(t;x,\psi)$ the set of $p\in U(t)$ that realize the maximum
of $\ct{H}_\gamma(x,t,p,\lambda,\psi)$. For all
$\gamma,\lambda\in\mm{T},(x,\psi)\in\mm{X}\times\mm{X}$,
 the compact-valued map $t\mapsto\ct{P}_{\gamma,\lambda}(t;x,\psi)$
 has a measurable selector by virtue
of~\cite[Theorem 3.7]{select}. Then, by \cite[Lemm 2.3.11]{tovst1},
for an arbitrary function~$(x,\psi)\in C(\mm{T},\mm{X}\times\mm{X})$
 the map $t\mapsto \ct{P}_{\gamma,\lambda}(t;(x,\psi)(t))$  
also  has a measurable selector.
Therefore, since relation~\rref{sys_max_} also depends on $x,\psi$
and on the parameters~$\gamma$ and~$\lambda$ upper semicontinuously, and all the relations are integrally bounded on bounded sets,
 by virtue of~\cite[Theorem 3.5.6]{tovst1}, on each finite interval
 for the funnels of solutions of~\rref{sys_x}--\rref{sys_psi},\rref{sys_max_}, we have upper semicontinuity by
 $\gamma,\lambda$.
  In particular, since $\gamma_n\to 0$ and $\lambda^n\to \lambda^0$,
  the upper limit of the compacts $\fr{Z}^n$
   is included in $\td{\fr{Z}}$. Hence,
   $(x^{00},\eta^0,\lambda^0,\psi^0)\in\td{\fr{Z}}$.

 On the other side, by $w\in(Null)(u^0)$ and by optimality of~$\eta^n$, $u^0$ for their problems, we obtain
\begin{equation*}
 \td{J}^{\eta^n}(\tau_n)-\gamma_n \fr{L}_w[\eta^n](\tau_n)\geq
 J^{u^0}(\tau_n)\geqtref{555} \td{J}^{\eta^n}(\tau_n)-\gamma_n^2
\end{equation*}
therefore, we have
   $\displaystyle
      \gamma_n \fr{L}_w[\eta^n](\tau_n)
 \leq \gamma_n^2.$
 By
    $\td{u}^0|_{[\tau_n,\infty\ra}=
   \eta^n|_{[\tau_n,\infty\ra},$
we obtain
\begin{equation}
   \label{to_w}
   \fr{L}_w[\eta^n](t)\leq \gamma_n
   \qquad\forall t\in\mm{T}.
\end{equation}
For each $t\in\mm{T}$, passing to the limit as $n\to\infty$, we
obtain that $\fr{L}_w[\eta^0]\leq 0$; i.e., $\fr{L}_w[\eta^0](t)=0$
for all $t\in\mm{T}$. Since $w\in(Null)(u^0)$, we have
$\eta^0=\td{u}^0$ a.e. on $\mm{T}$, hence
   $x^{00}\equiv x^0$ and $(x^0,u^0,\lambda^0,\psi^0)\in{\fr{Z}}$.
    Moreover, from~\rref{to_w}, we have $||\fr{L}_w[\eta^n]||_C\to 0$.
    \bo
\end{subequations}
\begin{acknowledgements}
I am grateful to an anonymous referee for helpful comments. I
would like to express my gratitude to  S.~M.Aseev, A.~G.~Chentsov,
A.~M.~Tarasyev, N.~Yu.~Lukoyanov, and Yu.~V.~Averboukh for valuable
discussion in course of writing this article. Special thanks to
Ya.~V.~Salii for the translation.

This supported by the Russian Foundation for Basic Research (RFBR) under
grant No~12-01-00537.
\end{acknowledgements}

\bibliographystyle{spmpsci}      


\end{document}